\documentclass[12pt]{article}

\usepackage[english]{babel}
\usepackage{graphicx}
\usepackage{enumerate}
\usepackage{amsfonts}
\usepackage{amssymb}
\usepackage[centertags]{amsmath}
\usepackage{amsthm}
\usepackage[usenames]{color}

\def\squarebox#1{\hbox to #1{\hfill\vbox to #1{\vfill}}}
\renewcommand{\qed}{\hspace*{\fill}
\vbox{\hrule\hbox{\vrule\squarebox{.667em}\vrule}\hrule}\smallskip}
\newtheorem{teorema}{Theorem}[section]
\newtheorem{lema}[teorema]{Lemma}
\newtheorem{corolario}[teorema]{Corollary}
\newtheorem{proposicao}[teorema]{Proposition}
\newtheorem{definicao}[teorema]{Definition}
\newenvironment{demonstracao}{\noindent {\bf Proof:}}{\hfill $\qed $ \newline}
\newenvironment{remar}{\noindent {\bf Remark:}}{}
\newenvironment{exemplo}{\noindent {\bf Example:}}{}

\begin{document}

\title{Morse and Lyapunov  Spectra and Dynamics on Flag Bundles\thanks{
Research supported by FAPESP grant n$^{0}$ 02/10246-2}}
\author{Luiz A. B. San Martin\thanks{%
Supported by CNPq grant n$^{\protect\underline{\circ }}$ 305513/03-6} \and %
Lucas Seco\thanks{%
Supported by FAPESP grant n$^{\protect\underline{\circ }}$ 04/00392-7} \\
%EndAName
Instituto de Matem\'{a}tica\\
Universidade Estadual de Campinas\\
Cx.Postal 6065, 13.081-970 Campinas-SP, Brasil}
\maketitle

\begin{abstract}
This paper studies characteristic exponents of flows  in relation
with the dynamics of flows on flag bundles. The starting point is
a flow on a principal bundle with semi-simple group $G$.
Projection against the Iwasawa decomposition $G = KAN$ defines an
additive  cocycle over the flow  with values in $\frak{a} = \log A$.
Its Lyapunov exponents (limits along trajectories) and Morse exponents
(limits along chains)  are studied. It is proved a
symmetric property of these  spectral sets, namely invariance
under the  Weyl group. It is proved also that these sets are
located in certain Weyl chambers,  defined from the dynamics
on the associated flag bundles. As a special case linear flows on
vector bundles are considered.
\end{abstract}



\section{Introduction}

In this paper we study characteristic exponents (Morse and Lyapunov
vector spectra) associated to continuous flows evolving on principal
bundles. The purpose is to relate these exponents to the dynamics on
the flag bundles. We work with invariant flows on a general
principal bundle whose structure group $G$ is semi-simple or
reductive. The objects to be considered are defined intrinsically
from $G$ and the principal bundle. Thus our characteristic exponents
take values on a vector subspace $\frak{a}$ of the Lie algebra
$\frak{g}$ of $G$, while the flag bundles in question have fiber
$\Bbb{F}$, a compact homeneous space of $G$ (generalized flag
manifolds of $G$). When $G$ is a linear group $\frak{a}$ becomes a
subspace of diagonal matrices and $\Bbb{F}$ a manifold of flags of
subspaces.

The vector valued spectra are defined via the cocycle over the
flag bundle obtained by projection against an Iwasawa
decomposition of $G$. They measure the exponential growth ratio of
the flow the same way as the usual one-dimensional spectra
measures the growth ratio of linear flows on vector bundles.
(Actually the one-dimensional spectrum can be recovered from the
vector valued one, via a representation of $G$, see Section
\ref{secvecbundl}.) This is the set up of  Kaimanovich \cite{ka}
to the proof  of the multiplicative ergodic theorem of Oseledets.

The Morse spectrum set is a concept of growth ratio along chains
introduced by Colonius-Kliemann \cite{ck1}, \cite{ck} (see also
Colonius-Fabri-Johnson \cite{cfj} for a vector valued version).
Each one of these sets depend on a chain component of a flow. In
 our case the Morse spectra are defined over chain components of the flow
  on an associated flag bundle.

Our purpose here is to relate the geometry of the vector Lyapunov
and Morse spectra with the geometry of the chain components on the
flag bundle, which were studied in \cite{smbflow} and \cite{msm}.

We better explain our results with the concrete example where
$G=\mathrm{Sl}\left( 3,\Bbb{R}\right) $. Let $X$ be a compact
Hausdorff space endowed with a continuous flow $\left( t,x\right)
\mapsto t\cdot x$ ($t\in \Bbb{T}=\Bbb{Z} $ or $\Bbb{R}$) and
suppose that $\phi _{t}\left( x,g\right) =\left( t\cdot x,\rho
\left( t,x\right) g\right) $ is a flow on $Q=X\times G$ with $\rho
$ a continuous cocycle taking values in $G$. Let $A\subset G$ be
the subgroup of diagonal matrices (positive entries) and $N$ the
subgroup of upper triangular unipotent matrices. Then the Iwasawa
decomposition of $G$ reads  $G=KAN$ with $K=\mathrm{SO}\left(
3\right) $.  For $x\in X$ and $u\in K$ we write $\rho \left(
t,x\right) u=k_{t}\left( x,u\right) a_{t}\left( x,u\right)
n_{t}\left( x,u\right) \in KAN$ and look at the limit behavior as
$t\rightarrow \infty $ of the matrix $\mathsf{a}\left(
t,x,u\right)/t =\log a_{t}\left( x,u\right)/t $ in the subspace
$\frak{a} = \log A$.

This is done in combination with the asymptotics of the compact
part  $k_{t}\left( x,u\right) $. Via the action of $G$ on $K\simeq
G/AN$ the
component $k_{t}\left( x,u\right) $ becomes a flow on $K$ and $\mathsf{a}%
\left( t,x,k\right) $ an additive cocycle over this flow. Thus the
characteristic exponents defined by $\mathsf{a}$ depend on $\left(
x,k\right) \in X\times K$. Nevertheless it is more convenient to avoid
ambiguities and factor out the subgroup $M\subset K$ of diagonal matrices
with $\pm 1$ entries and get $\mathsf{a}\left( t,x,b\right) $ as a cocycle
over the flow $\phi ^{\Bbb{F}}$ induced on $X\times \Bbb{F}$, where $\Bbb{F}%
\simeq K/M\simeq G/MAN$ is the manifold of complete flags in
$\Bbb{R}^{3}$.

Now let $\mathcal{M}$ be a chain component of the flow $\phi ^{\Bbb{F}}$ on $%
X\times \Bbb{F}$. Its associated Morse spectrum set $\Lambda _{\mathrm{Mo}%
}\left( \mathcal{M}\right) \subset \frak{a}$ is defined by evaluating the
cocycle $\mathsf{a}\left( t,x,b\right) $ on chains in $\mathcal{M}$ taking
into account the jumps of the chains (see Definition \ref{defmorseexpon}
below). The Lyapunov spectrum set $\Lambda _{\mathrm{Ly}}\left( \mathcal{M}%
\right) $ is the set of limits $\lim \mathsf{a}\left(
t,x,b\right)/t$
along the trajectories in $\mathcal{M}$. By general facts $\Lambda _{\mathrm{%
Mo}}\left( \mathcal{M}\right) $ is a compact convex set which contains $%
\Lambda _{\mathrm{Ly}}\left( \mathcal{M}\right) $, which in turn
contains the set of extremal points of $\Lambda _{\mathrm{Mo}}\left(
\mathcal{M}\right) $ (see Section \ref{secmorgeral} below).

The chain components in $X\times \Bbb{F}$ were described in
\cite{smbflow}, \cite{msm} with the assumption that the flow on
$X$ is chain transitive. The picture is the following: There
exists $H_{\phi }\in \frak{a}$ and a continuous map $x\in X\mapsto
f_{\phi }\left( x\right) =g_{x}H_{\phi }g_{x}^{-1}$, $g_{x}\in G$,
into the set of  conjugates of $H_{\phi }$ such that each chain
component $\mathcal{M}$ has the form $\bigcup_{x\in X}\{x\}\times \mathrm{fix%
}\left( f_{\phi }\left( x\right) \right) $. Here
$\mathrm{fix}\left( f_{\phi }\left( x\right) \right) $ is a
connected component of the fixed point set of the action on
$\Bbb{F}$ of  $\exp f_{\phi }\left( x\right) $ (the connected
components can be labelled accordingly so that they match each
other). For example if the diagonal matrix $H_{\phi }$ has
distinct eigenvalues $\lambda _{1}>\lambda _{2}>\lambda _{3}$ with
eigenvectors  $e_{1},e_{2},e_{3}$ then there are $3!=6$ isolated
fixed points, namely the flag
$b_{0}=(\mathrm{span}\{e_{1}\}\subset
\mathrm{span}\{e_{1},e_{2}\})$ together with those obtained from
$b_{0}$ by permuting the basic vectors. A similar description
holds for $f_{\phi }\left( x\right) $, $x\in X$, allowing to label
the fixed points, as well as the chain components, with the
permutation group on three letters $S_3$, which is the Weyl group
$\mathcal{W}$ for $\mathrm{Sl}\left( 3,\Bbb{R}\right) $.
In this case there are $6$ distinct chain
components denoted by $\mathcal{M}_{w}$ with $w$ running over
$\mathcal{W} = S_3$.

The chain components are also labelled by the permutation group
even if  $H_{\phi }$ has repeated eigenvalues, but now some of the
$\mathcal{M}_{w}$ can merge. In our $\mathrm{Sl}\left(
3,\Bbb{R}\right) $ example there are four possibilities for
$H_{\phi }$, namely

\begin{enumerate}[(a)]
\item  $H_{\phi }=\mathrm{diag}(\lambda _{1}>\lambda _{2}>\lambda _{3})$
with six chain components $\mathcal{M}_{w}$, $w\in \mathcal{W}$.

\item  $H_{\phi }=\mathrm{diag}(\lambda _{1}=\lambda _{2}>\lambda _{3})$
with three chain components $\mathcal{M}_{1}=\mathcal{M}_{\left(
12\right) }$,
 $\mathcal{M}_{\left( 23\right) }=\mathcal{M}_{\left( 123\right) }$ and
 $ \mathcal{M}_{\left( 132\right) }=\mathcal{M}_{\left( 13\right) }$.

\item  $H_{\phi }=\mathrm{diag}(\lambda _{1}>\lambda _{2}=\lambda _{3})$
with three chain components $\mathcal{M}_{1}=\mathcal{M}_{\left( 23\right) }$%
, $\mathcal{M}_{\left( 12\right) }=\mathcal{M}_{\left( 132\right) }$ and $%
\mathcal{M}_{\left( 123\right) }=\mathcal{M}_{\left( 13\right) }$.

\item  $H_{\phi }=\mathrm{diag}(\lambda _{1}=\lambda _{2}=\lambda _{3}) = 0$ and
the flow is chain transitive on $X \times {\mathbb F}$.
\end{enumerate}

\noindent (In all the cases $\mathcal{M}_{1}$ is the only attractor
component and $\mathcal{M}_{\left( 13\right) }$ the only repeller. The
equalities in the second and third cases are due to the cosets $\mathcal{W}%
_{\phi }w$ where $\mathcal{W}_{\phi }$ is the subgroup of
$\mathcal{W}$ fixing $H_{\phi }$.)

We say that the matrix $H_{\phi }$ is the block form of $\phi $. It
is determined by the dynamics on the flag bundles and is a
combination of the parabolic type of the flow and of the reversed
flow (see \cite{smbflow}, \cite{msm} and Section \ref{secpartydyn}
of the present article).

The results of this paper provide a similar picture for the Morse
and Lyapunov spectra over the chain components in terms of the Weyl
group, the Weyl chambers and the block form $H_{\phi }$ of $\phi $.

For $G=\mathrm{Sl}\left( 3,\Bbb{R}%
\right) $ the roots of $\frak{a}$ are the functionals $\alpha
_{ij}\left( \mathrm{diag}(a_{1},a_{2},a_{3})\right) =a_{i}-a_{j}$,
for $i\neq j \in \{1,2,3\}$. Their kernels are three straight
lines containing the nonregular matrices in ${\mathfrak a}$, i.e.,
those having repeated eigenvalues.  The complement to the lines is
the union of $3!=6$ Weyl chambers where the Weyl group acts simply
and transitively  (see Figure \ref{fig:espectro-possib-1}). A
chamber
is determined by inequalities $a_{i}>a_{j}>a_{k}$ and we denote it by $%
\mathcal{C}_{w}$ if $w\in \mathcal{W}$ is the permutation sending
$\left( 123\right) $ to $\left( ijk\right) $.

\begin{figure}[t]
    \begin{center}
        \includegraphics{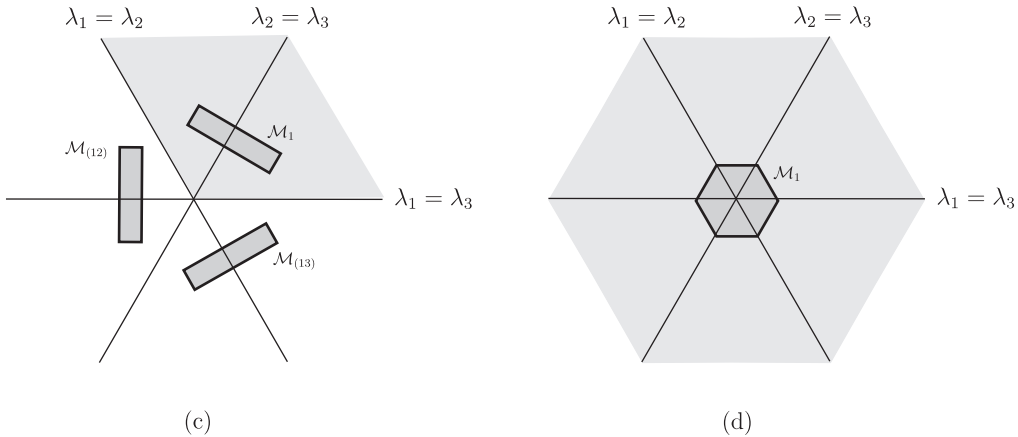}
    \end{center}
    \caption{\label{fig:espectro-possib-2}
    Morse spectra for flows with
    $G = \mathrm{Sl}(3,{\mathbb R})$: block forms (c) and (d).}
\end{figure}

The main results of this paper prove the following estimates and
symmetry properties of the vector valued Morse and Lyapunov
spectra (see Figures \ref{fig:espectro-possib-1} and
\ref{fig:espectro-possib-2}):

\begin{figure}[t]
    \begin{center}
        \includegraphics{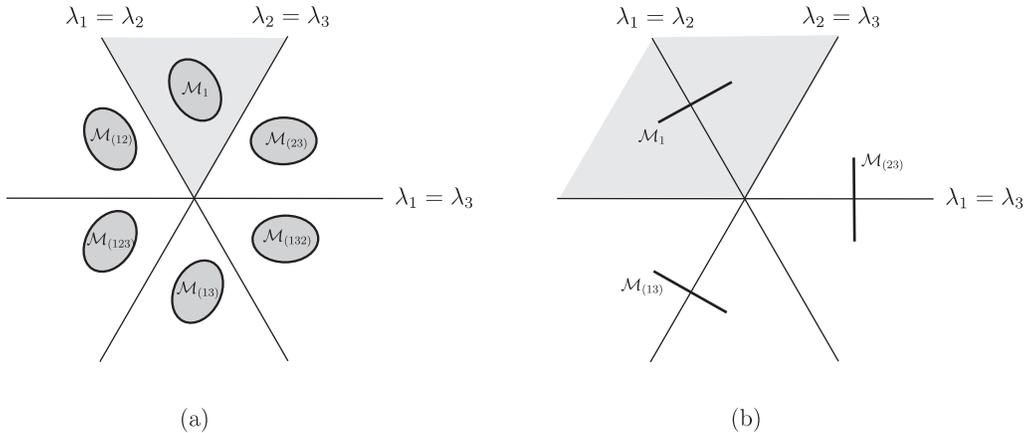}
    \end{center}
    \caption{\label{fig:espectro-possib-1}
    Morse spectra for flows with
    $G = \mathrm{Sl}(3,{\mathbb R})$: block forms (a) and (b), each spectrum is
labeled by its corresponding component.}
\end{figure}

\begin{enumerate}
\item  The union $\bigcup_{w\in \mathcal{W}}\Lambda _{\mathrm{Mo}}\left(
\mathcal{M}_{w}\right) $ of the several Morse spectra is invariant under the
Weyl group. The same holds for the Lyapunov spectra $\bigcup_{w\in \mathcal{W}%
}\Lambda _{\mathrm{Ly}}\left( \mathcal{M}_{w}\right) $.

\item  The spectra $\Lambda _{\mathrm{Mo}}\left( \mathcal{M}_{1}\right) $
and $\Lambda _{\mathrm{Ly}}\left( \mathcal{M}_{1}\right) $ of the attractor
component are invariant under the subgroup $\mathcal{W}_{\phi }$ fixing $%
H_{\phi }$.

\item  $\Lambda _{\mathrm{Mo}}\left( \mathcal{M}_{1}\right) $ is contained
in the interior of $\bigcup_{w\in \mathcal{W}_{\phi }}\mathrm{cl}\mathcal{C}%
_{w}$. This means that the multiplicities of the eigenvalues of any
Morse exponent in $\Lambda _{\mathrm{Mo}}\left(
\mathcal{M}_{1}\right) $ do not
exceed those of $H_{\phi }$. The same estimate holds for $\Lambda _{\mathrm{%
Ly}}\left( \mathcal{M}_{1}\right) $, since $\Lambda
_{\mathrm{Mo}}\left(
\mathcal{M}_{w}\right) $ is the convex closure of $\Lambda _{\mathrm{Ly}%
}\left( \mathcal{M}_{w}\right) $.

\item  $\Lambda _{\mathrm{Mo}}\left( \mathcal{M}_{1}\right) $ intercepts the
closure of every chamber $\mathcal{C}_{w}$, $w\in \mathcal{W}_{\phi }$.
Hence there exists a Morse exponent in $\Lambda _{\mathrm{Mo}}\left(
\mathcal{M}_{1}\right) $ whose eigenvalues have the same pattern of
multiplicities as $H_{\phi }$. $\Lambda _{\mathrm{Ly}}\left( \mathcal{M}%
_{1}\right) $ also meets these chambers, but since it is not
necessarily convex, it may happen that every Lyapunov exponent has
less multiplicities than $H_{\phi }$.

\item  $\Lambda _{\mathrm{Mo}}\left( \mathcal{M}_{w}\right) =w^{-1}\Lambda _{%
\mathrm{Mo}}\left( \mathcal{M}_{1}\right) $ and $\Lambda _{\mathrm{Ly}%
}\left( \mathcal{M}_{w}\right) =w^{-1}\Lambda _{\mathrm{Ly}}\left( \mathcal{M%
}_{1}\right) $. Combining this with the previous statement it follows that
distinct Morse spectra do not overlap. (This fact is not true for linear
flows on vector bundles as shown in \cite{ck}, Example 5.5.11. The point
here is that the vector bundle Morse spectra are images under linear maps of
our spectra, see Section \ref{secvecbundl}. Overlappings of the images may occur.)
\end{enumerate}

The last statement says that the whole Morse spectra can be read off
from the spectrum of the attractor component. This phenomenon is
already present in the analysis of the chain components in flag
bundles, whose main properties are governed by the extremal
(attractor and repeller) components.

These results show an intimate relation between the spectra and the
dynamics on the flag bundles. From the chain components on the
bundles we get the symmetries of the spectra. Conversely, from the
Morse spectra (of the attractor component only) we can recover the
block form $H_{\phi }$ of the flow, and hence its parabolic type.
This in turn tells the number of chain components, their geometry,
their Conley indices, etc.

The underlying structure behind this relationship is a ``block
decomposition'' of the flow in the following sense: From the chain
components on the flag bundles one can build a $\phi
_{t}$-invariant subbundle $Q_{\phi }\subset Q$ whose structure
group is the centralizer  $Z_{H_{\phi }}$ of the block form
$H_{\phi }$ (see Section \ref{secpartydyn} below). If
$G=\mathrm{Sl}\left( n,\Bbb{R}\right) $ then $Z_{H_{\phi }}$ is a
subgroup of block diagonal matrices, hence we say that $Q_{\phi }$
is the block reduction of $\phi $. The above results say that the
spectra have the same block structure as the block reduction, as
expected. In case the bundles $Q$ and $Q_{\phi }$ are trivial then
the reduction amounts a cohomology of cocycles reducing the
original cocycle to one taking values in the subgroup
$Z_{H_{\phi}}$. Hence this reduction points towards a Jordan
decomposition of the flow, in the sense of Arnold-Cong-Oseledets
\cite{aco}.

\begin{figure}[t]
    \begin{center}
        \includegraphics{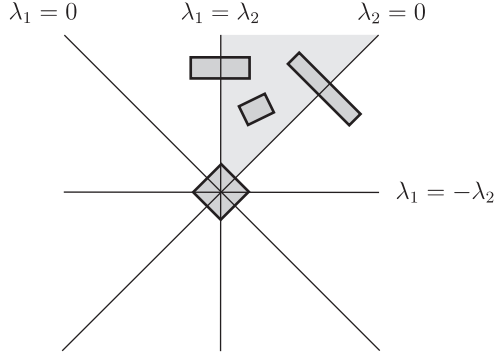}
    \end{center}
    \caption{\label{fig:espectro-sp-4}
    For flows with $G = \mathrm{Sp}(4,{\mathbb R})$ the spectrum lives here.}
\end{figure}

A similar picture to the one described in the above examples holds
for invariant flows on a general principal bundle whose structure
group $G$ is a noncompact semi-simple Lie group with finite center
(the Morse spectra for the  symplectic group $\mathrm{Sp}\left(
4,{\mathbb R} \right) $ are depicted in  Figure
\ref{fig:espectro-sp-4} where ${\mathfrak a}$ consists of the
diagonal matrices ${\rm diag}(\lambda_1, \lambda_2, -\lambda_1,
-\lambda_2)$). We develop the setup and state our results in the
generality of semi-simple Lie groups. These results can be
extended to reductive groups (e.g.\ $\mathrm{Gl}\left(
n,\Bbb{R}\right) $). But now one must add a central component to
$\frak{a}$ and hence to $\Lambda _{\mathrm{Mo}}$. There is not
much to say about this central component unless that it is the
same for every chain component $\mathcal{M}$.

\section{Preliminary notation\label{sec:prelim}}

For the theory of semi-simple Lie groups and their flag manifolds
we refer to Duistermat-Kolk-Varadarajan \cite{dkv}, Helgason
\cite{Helgason} and Warner \cite{w}. To set notation let $G$ be a
noncompact  semi-simple Lie group  with  Lie algebra $\frak{g}$.
We assume throughout that $G$ has finite center. Fix a Cartan
involution $\theta $ of $\frak{g}$ with Cartan decomposition
$\frak{g}=\frak{k}\oplus \frak{s}$. The form   $B_{\theta }\left(
X,Y\right) =-\langle X,\theta Y\rangle $, where $\langle \cdot
,\cdot \rangle $ is the Cartan-Killing form of $\frak{g}$, is an
inner product.

Choose   a maximal abelian subspace $\frak{a}\subset \frak{s}$ and
a Weyl chamber $\frak{a}^{+}\subset \frak{a}$. We  let $\Pi $ be
the set of roots of $\frak{a}$,  $\Pi ^{+}$ the positive roots
corresponding to $\frak{a}^{+}$ and $\Sigma $ the set of simple
roots in $\Pi ^{+}$.  The associated Weyl group is
$\mathcal{W}=M^{*}/M$ where $M^{*}$ and $M$ are  the normalizer
and the centralizer of $A$, respectively.  We write $\frak{m}$ for
the  Lie algebra of $M$.

The Iwasawa decomposition reads  $\frak{g}=\frak{k}\oplus
\frak{a}\oplus \frak{n}^{+}$ with  $\frak{n}^{+}=\sum_{\alpha \in
\Pi ^{+}}\frak{g}_{\alpha }$  where $\frak{g}_{\alpha }$ is the
root space associated to $\alpha $. Put $\frak{n}^{-}=\sum_{\alpha
\in -\Pi ^{+}}\frak{g}_{\alpha }$. As to  the global
decompositions we write  $G=KS$ and  $G=KAN$ with $K=\exp
\frak{k}$, $S=\exp \frak{s}$, $A=\exp \frak{a}$ and $N=\exp
\frak{n}^{+}$. The  standard minimal parabolic subalgebra is
defined by
\[
{\frak p}={\frak m}\oplus {\frak a}\oplus {\frak n}^{+}
\]

Associated to  a subset of simple roots $\Theta \subset \Sigma $
there are several Lie  algebras and groups (cf.  \cite{w}, Section
1.2.4): We write $\frak{g}\left( \Theta \right) $ for the
(semi-simple)  Lie subalgebra
generated by $\frak{g}_{\alpha }$, $\alpha \in \Theta $, and put $\frak{a}%
\left( \Theta \right) =\frak{g}\left( \Theta \right) \cap \frak{a}$ and $%
\frak{n}^{\pm }\left( \Theta \right) =\frak{g}\left( \Theta
\right) \cap \frak{n}^{\pm }$. $G\left( \Theta \right) $ is the
connected group with Lie algebra $\frak{g}\left( \Theta \right) $
and $A\left( \Theta \right) =\exp
\frak{a}\left( \Theta \right) $, $N^{\pm }\left( \Theta \right) =\exp \frak{n%
}^{\pm }\left( \Theta \right) $. Also, $\frak{a}_{\Theta }=\{H\in
\frak{a}:\alpha \left( H\right) =0$, $\alpha \in \Theta \}$ and
$A_{\Theta }=\exp \frak{a}_{\Theta }$.

We let  $Z_{\Theta }$ be  the  centralizer of $\frak{a}_{\Theta }$
in $G$ and $K_{\Theta }=Z_{\Theta }\cap K$.  It decomposes as
$Z_{\Theta }=MG(\Theta )A_{\Theta }$ which implies that $Z_{\Theta
}=K_{\Theta }(S\cap Z_{\Theta })$ is a Cartan decomposition  while
$Z_{\Theta }=K_{\Theta }AN(\Theta )$ is an Iwasawa decomposition
of $Z_{\Theta }$ (which is a reductive Lie group).

The standard parabolic subagebra $\frak{p}_{\Theta
}=\frak{n}^{-}\left( \Theta \right) \oplus \frak{p}$ and the
corresponding standard parabolic subgroup $P_{\Theta }$ is the
normalizer of $\frak{p}_{\Theta }$ in $G$. The flag manifold  of
type $\Theta $ is $\Bbb{F}_{\Theta }=G/P_{\Theta }$, which
identifies with the set of conjugates of $\frak{p}_{\Theta }$. The
empty set $\Theta =\emptyset $ gives the minimal parabolic
subgroup $P=P_{\emptyset }$ , which is $P=MAN$.

We write $\frak{p}_{\Theta }^{-}=\theta (\frak{p}_{\Theta })$ for
the parabolic subalgebra opposed to $\frak{p}_{\Theta }$. It is
conjugate to the parabolic subalgebra $\frak{p}_{\Theta ^{*}}$
where $\Theta ^{*}=-(w_{0})\Theta $ and $w_{0}$ is the principal
involution of $\mathcal{W} $ (the one that takes $\Sigma $ to
$-\Sigma $). More precisely, $\frak{p}_{\Theta
}^{-}=\overline{w}_{0}\frak{p}_{\Theta ^{*}}$ if
$\overline{w}_{0}\in M^{*}$ is a representative of $w_{0}$. If
$P_{\Theta }^{-}$ is the parabolic subgroup associated to
$\frak{p}_{\Theta }^{-}$ then
\[
Z_{\Theta }=P_{\Theta }\cap P_{\Theta }^{-},
\]
and $P_{\Theta }^{-}=N_{\Theta }^{-}Z_{\Theta }$.

The subset $\Theta $ singles out the subgroup of the Weyl group $\mathcal{W}%
_{\Theta }=\left( \mathcal{W}\cap P_{\Theta }\right) /M=\left( \mathcal{W}%
\cap Z_{\Theta }\right) /M$. Alternatively $\mathcal{W}_{\Theta }$
is  the subgroup generated by the reflections with respect to the
roots $\alpha \in \Theta $.

For $H\in \frak{a}$ let  $\Theta _{H}=\{\alpha \in \Sigma :\alpha
(H)=0\}$. Then we denote  the above subalgebras an groups with $H$
instead of $\Theta _{H}$, for example,
$\frak{p}_{H}=\frak{p}_{\Theta _{H}}$, etc. (cf. Section
\ref{secpartydyn}). In this case $Z_{H}=Z_{\Theta _{H}}$ and $%
K_{H}=K_{\Theta _{H}}$ are  the centralizer of $H$ in $G$ and $K$
respectively. Also,  $\mathcal{W}_{H}=\mathcal{W}_{\Theta _{H}}$
is the subgroup of $\mathcal{W}$ fixing $H$. Conversely, for
appropriately chosen elements $H\in \frak{a}$ such that $\Theta
_{H}=\Theta $, the parabolic subgroups, flags, etc. can be written
with reference to $H$.

\section{Morse spectrum of a vector valued cocycle\label{secmorgeral}}

In this section we discuss the abstract concepts of vector valued
cocycles and their Lyapunov and Morse characteristic exponents. We
mostly recall the
results of Colonius-Kliemann \cite{ck} and Colonius-Fabbri-Johnson \cite{cfj}%
, which are stated in the continuous-time setting for flows on
metric spaces. We indicate here the mild modifications needed for
discrete-time flows on more general topological spaces. Detailed
proofs will be given elsewhere \cite{sl}.

Let $E$ be a compact Hausdorff space and $\phi :\Bbb{T}\times E\to E$ a
continuous flow with $\Bbb{T}=\Bbb{R}$ or $\Bbb{Z}$. We denote $\phi
(t,x)=\phi _{t}(x)=t\cdot x$. If $V$ is a finite dimensional normed vector
space, a $V$-valued cocycle over $E$ is a continuous map $a:\Bbb{T}\times
E\to V$ with
\[
a(t+s,x)=a(t,s\cdot x)+a(s,x).
\]
Two $V$-valued cocycles $a$ and $b$ over $E$ are cohomologous with
cohomology $h:E\to V$ if $a(t,x)+h(x)=h(t\cdot x)+b(t,x)$, where $h$
is a continuous map. A $V$-valued cocycle $a(t,x)$ over $E$ defines
a flow on $E\times V$ by
\[
t\cdot (x,v)=(t\cdot x,\,a(t,x)+v).
\]

The product $E\times V$ can be viewed as principal $V$-bundle over $E$,
where the aditive group $V$ acts on the right and leaves invariant the above
flow, which is thus made of automorphisms of the bundle. Conversely, let $%
\pi :P\to E$ be a principal bundle with vector structure group $V$ and $\phi
_{t}$ a flow of automorphisms of $P$. It is well known that $P$ is a trivial
bundle (see e.g. Kobayashi-Nomizu \cite{kn}). To each global continuous
section $\chi :E\to P$ there corresponds a trivialization of $P$ and a
continuous $V$-cocycle $a_{\chi }(t,x)$ over $E$ given by
\[
t\cdot \chi (x)=\chi (t\cdot x)\cdot a_{\chi }(t,x).
\]
We note that any two cocycles coming this way from sections are
cohomologous. In fact, if $\eta :E\to P$ is another section of $P$ and $%
a_{\eta }$ the corresponding cocycle, then there exists a continuous map $%
h:X\to V$ satisfying $\chi (x)=\eta (x)\cdot h(x)$, and it follows that
\[
a_{\eta }(t,x)+h(x)=h(t\cdot x)+a_{\chi }(t,x),
\]
so that the map $h$ realizes a cohomology between $a_{\eta }$ and
$a_{\chi } $. Vector cocycles coming from such trivializations of
principal bundles will show up below.

Given $x\in E$, $T\in \Bbb{T}$, the finite-time Lyapunov exponent of the
cocycle $a$ at $(x,T)$ is
\[
\lambda _{T}(x)=\frac{1}{T}a(T,x),
\]
while the Lyapunov exponent of \emph{$a$ }at (in the direction of) \emph{$x$}
is
\[
\lambda (x)=\lim_{T\rightarrow +\infty }\lambda _{T}(x),
\]
if the limit exists. The Lyapunov spectrum of a subset\emph{\
$Y\subset E$} is defined by
\[
\Lambda _{\mathrm{Ly}}(Y)=\{\lambda (y)\in V:\,y\in Y\,\mathrm{and}\,\text{$%
\lambda (y)\,$}\mathrm{exists}\}.
\]

The exponential growth ratio along chains was introduced by
Colonius-Kliemann \cite{ck1}, who called them Morse exponents. We
use the approach to chains for flows in topological spaces as
developed in Patr\~{a}o \cite{patr}. Thus let $\mathcal{O}$ be a
family of open coverings of the compact Hausdorff space $E$, which
is admissible in the sense of \cite {patr}. If $\mathcal{U}\in
\mathcal{O}$ is given then two points $x,y\in E$ are said to
$\mathcal{U}$-close if there exists an open set $A\in \mathcal{U}$
such that $x,y\in A$.

A $(\mathcal{U},T)$-chain ($\mathcal{U}\in \mathcal{O}$, $T\in \Bbb{T}$)
$\zeta $ from $x$ to $y$ in $E$ is a finite sequence of points $%
x=x_{0},x_{1},\ldots ,x_{N}$ in $E$ and times $t_{0},\,t_{1},\ldots ,t_{N-1}
$ in $\Bbb{T}$ with $t_{i}\geq T$ such that $t_{i}\cdot x_{i}$ and $%
\,x_{i+1} $ are $\mathcal{U}$-close for $i=0,\ldots ,N-1$. We say that the
chain $\zeta $ belongs to the subset $\mathcal{M}\subset E$ if the initial
and end points $x,y$ belong to $\mathcal{M}$ (but not necessarily the
intermediate points).

Let $a$ as above be a cocycle over $E$.

\begin{definicao}
\label{defmorseexpon}The finite time Morse exponent $\lambda (\zeta )$ of
the $(\mathcal{U},T)$-chain $\zeta $ is
\[
\lambda (\zeta )=\frac{1}{T(\zeta )}\sum_{j=0}^{N-1}a(T_{j},x_{j}),
\]
where $T(\zeta )=\sum_{j=0}^{N-1}T_{j}$ is the total time of $\zeta $. We
denote by $\Lambda _{\mathrm{Mo}}(\mathcal{M};{\mathcal{U}},T)$ the set of $%
\left( \mathcal{U},T\right) $-chains belonging to $\mathcal{M}$ and define
the Morse spectrum of $\mathcal{M}$ to be
\[
\Lambda _{\mathrm{Mo}}(\mathcal{M})=\bigcap \{\mathrm{cl}\Lambda _{\mathrm{Mo%
}}(\mathcal{M};{\mathcal{U}},T):\,\mathcal{U}\in \mathcal{O},\,T>0\}.
\]
\end{definicao}

Note that the family of coverings $\mathcal{O}$ is an ingredient in the
definition of $\Lambda _{\mathrm{Mo}}(\mathcal{M})$. However this set does
not depend on the specific choice of $\mathcal{O}$ (see Corollary \ref
{teo-espectro-independe-da-cob} below).

It follows immediately from the definitions that
\[
\lambda (\zeta )=\sum_{j=0}^{N-1}\frac{T_{j}}{T(\zeta )}\lambda
_{T_{j}}(x_{j}),
\]
so that each finite-time Morse exponent is a convex combination of
finite-time Lyapunov exponents. That an analogous result is valid for the
limit Morse exponents is one of the central results on the structure of
Morse exponents which is stated next, together with other results.

\begin{teorema}
\label{propespectromorse}

Let $\mathcal{M}\subset E$ be a chain component (that is, a maximal
chain transitive set).

\begin{enumerate}
\item  Two cohomologous cocycles have the same Lyapunov and Morse spectra.

\item The Morse spectrum $\Lambda _{\mathrm{Mo}}(\mathcal{M})$ of a continuous time flow
 coincides with any of its discretezations (restriction to $c{\mathbb Z}$,
 $c\in  {\mathbb R}$).

\item  Suppose that $\mathcal{M}$ contains the $\omega $-limit set $%
\omega (x)$ of $x\in E$. Then $\lambda (x)\in \Lambda _{\mathrm{Mo}}(%
\mathcal{M})$, if $\lambda (x)$ exists.

\item  $\Lambda _{\mathrm{Mo}}(\mathcal{M})$ is a
nonempty compact convex subset of $V$.

\item  $\Lambda _{\mathrm{Mo}}(\mathcal{M})$ is the
set of integrals $\int q \, \mathrm{d}\Bbb{P}$ with $\Bbb{P}$
running through the probability measures supported in $\mathcal{M}$
which are invariant by the time-one flow $\phi_1$. Here $q:E\to V$
is the time-one map $q(x)=a(1,x)$.

\item  The extremal points of $\Lambda _{\mathrm{Mo}}(\mathcal{M})$ are Lyapunov
exponents so that $\Lambda _{\mathrm{Mo}}(\mathcal{M})$ is the
closed convex closure of $\Lambda _{\mathrm{Ly}}(\mathcal{M})$.

\end{enumerate}
\end{teorema}

\begin{demonstracao}
See \cite{cfj}.  For item (2) and proofs in the present more general context see \cite{sl}.
\end{demonstracao}

Some immediate applications of these results are given in what follows.

\begin{corolario}
\label{corol-espectro-lyap-morse}Let $\mathcal{M}_{j}$, $j=1,\ldots ,n$, be
the finest Morse decomposition of $\phi _{t}$ in $E$. Then
\[
\Lambda _{\mathrm{Ly}}(E)\subset \bigcup_{j=1}^{n}\Lambda _{\mathrm{Mo}}(%
\mathcal{M}_{j}).
\]
\end{corolario}

\begin{demonstracao}
In fact, the Morse components $\mathcal{M}_{j}$ are chain components
and contain all the $\omega $-limit sets in $E$.  The result then
follows from the above theorem.
\end{demonstracao}

\begin{corolario}
\label{teo-espectro-independe-da-cob}The Morse spectrum of a chain component
$\mathcal{M}$ does not depend on the family of coverings $\mathcal{O}$.
\end{corolario}

\begin{demonstracao}
In fact, by the above theorem we have $\Lambda _{\mathrm{Mo}}(\mathcal{M})=\mathrm{cl}\left( \mathrm{co}%
\Lambda _{\mathrm{Ly}}(\mathcal{M})\right) $.  Since $\Lambda _{\mathrm{Ly}}(%
\mathcal{M})$ does not depend on $\mathcal{O}$ the result follows.
\end{demonstracao}

\begin{corolario}
\label{propos:cociclos-relacionados}Let $\phi _{t}^{\prime }$ be a
flow on a compact Hausdorff space $E^{\prime }$ and let $b$ be a
$V^{\prime }$-vector cocycle over $E^{\prime }$. Suppose that there
is a map $\pi :E\to E^{\prime }$ and a linear map $p:V\to V^{\prime
}$ such that the cocycles $a$ and $b$ are related by
\[
b(t,\pi (x))=p(a(t,x)),\quad t\in \Bbb{T},\,x\in E,
\]
and that $\pi (\mathcal{M})$ is a chain transitive component in
$E^{\prime }$.  Then
\begin{equation}
\Lambda _{\mathrm{Mo}}(\pi (\mathcal{M}),b)=p(\Lambda _{\mathrm{Mo}}(%
\mathcal{M},a)),  \label{eq:cociclos-relacionados}
\end{equation}
where the notation indicates which cocycle is being considered.
\end{corolario}

\begin{demonstracao}
The Lyapunov exponents of $a$ and $b$ are easily seen to be related by (\ref
{eq:cociclos-relacionados}) with $\Lambda _{\mathrm{Ly}}$ in place of $%
\Lambda _{\mathrm{Mo}}$. The result then follows by the above
theorem.
\end{demonstracao}

\section{Decompositions of principal bundles and the $\frak{a}$-cocycle\label%
{secdecom}}

Let $Q\rightarrow X$ be a principal bundle whose structural group $G$ is a
reductive Lie group. In this section we build decompositions of $Q$ similar
to the Iwasawa and Cartan decompositions of $G$. We assume throughout that
the base space $X$ is paracompact. The right action of $G$ on $Q$ is denoted
by $q\mapsto q\cdot g$, $q\in Q$, $g\in G$.

\subsection{Iwasawa decomposition}

An Iwasawa decomposition $G=KAN$ of the structural group can be carried over
to a decomposition of the bundle, which we call the \textit{Iwasawa
decomposition} of $Q$ as well. The construction goes as follows. The group $AN$ is
diffeomorphic to an Euclidian space and since the base space
$X$ of $Q$ is paracompact, by general results on principal bundles there exists a
$K$-reduction of $Q$, that is, a subbundle $R\subset Q$ which is a
principal bundle with structural group $K$ (for a proof combine
Theorem I.5.7 with Proposition I.5.6 in Kobayashi-Nomizu
\cite{kn}). In case $G$ is a linear group such reduction can be
obtained from a Riemannian metric in an associated vector bundle.

For every $q\in Q$ there exists $g\in G$ and $r'\in R$ such that
$q=r'\cdot g$.  The $K$-component of the Iwasawa decomposition
$g=khn\in KAN$ of $g$ leaves $R$ invariant so that, by uniqueness of
the Iwasawa decomposition, it follows that every $q\in Q$ decomposes
uniquely as
\[
q=r\cdot hn, \quad r\in R,\quad hn\in AN,
\]
which exhibits the bundle $Q$ as a product $Q=R\cdot AN \approx R\times
A\times N$. We let
\[
\mathsf{R}:Q\to R,\quad q \mapsto r,\quad \qquad \mathsf{A}:Q\to
A,\quad q \mapsto h,
\]
be the associated projections. By the continuity of the Iwasawa
decomposition of $G$ and the local triviality of $Q$, it follows
that these projections are continuous. Moreover they satisfy the
following properties:

\begin{enumerate}
\item  $\mathsf{R}(r)=r$, $\mathsf{A}(r)=1$, when $r\in R$,

\item  $\mathsf{R}(q\cdot p)=\mathsf{R}(q)m$, $\mathsf{A}(q\cdot p)=\mathsf{A%
}(q)h$, when $q\in Q$, $p=mhn\in P=MA\!N$. In particular,
$\mathsf{A}(r\cdot p)=h$.
\end{enumerate}

In what follows we write for $q \in Q$,
\[
\mathsf{a}\left( q\right) =\log \mathsf{A}\left( q\right) \in
\frak{a}.
\]

We use the same notation for the Iwasawa decomposition of $g\in
G$, namely, $\mathsf{a}\left( g\right) =\log h $ if $g=uhn\in
KAN$. For future reference we note that by the second of the above
properties  we have
\begin{equation}\label{eq-propr-aditiva-iwasawa}
\mathsf{a}\left( q\cdot p\right) =\mathsf{a}\left( q\right)
+\mathsf{a}\left( p\right), \quad p \in P.
\end{equation}

Now we discuss the cocycle defined by $\mathsf{a}$. The first step is to get
actions on the $K$-subbundle $R$. Given an automorphism $\varphi \in \mathrm{%
Aut}(Q)$ define the map
\[
\varphi ^{R}:r\in R\mapsto \mathsf{R}(\varphi (r))\in R.
\]
We have
\[
\varphi ^{R}\circ \psi ^{R}=(\varphi \circ \psi )^{R}
\]
if $\psi \in \mathrm{Aut}(Q)$. In fact,  if $\psi \left( r\right)
=r_{1}\cdot g $, $r_{1}\in R$, $g\in AN$, then $\psi ^{R}\left(
r\right) =r_{1}$ so that by the above  properties
\[
(\varphi \circ \psi )^{R}\left( r\right) =\mathsf{R}(\varphi (\psi (r)))=%
\mathsf{R}(\varphi (r_{1})\cdot g)=\mathsf{R}(\varphi
(r_{1}))=\varphi ^{R}(r_{1})=\varphi ^{R}(\psi ^{R}(r)).
\]

We remark that $\varphi ^{R}$ is not a bundle morphism unless $R$
is  $\varphi $-invariant. Also, when $Q=G$ then $\varphi ^{R}$
reduces to the usual action of $G$ on $K$ through the Iwasawa
decomposition of $G$.

Let $\phi _{t}$ be the continuous flow of automorphisms of $Q$,
$t\in {\mathbb T} $. Then $\phi _{t}^{R}$ defines a continuous
flow in $R$. We use also the notation $\mathsf{a}$ for the cocycle
defined by the map $\mathsf{a}$, namely
\[
\mathsf{a}:\Bbb{T}\times R\to \frak{a},\quad
\mathsf{a}(t,r)=\mathsf{a}(\phi _{t}^{R} (r)),
\]
which is continuous by the continuity of $\phi $ and the projection $\mathsf{%
a}$. It follows from the properties of the Iwasawa decomposition of $Q$ that
$\mathsf{a}$ is a cocycle over $\phi _{s}^{R}$, that is,
\[
\mathsf{a}(t+s,r)=\mathsf{a}(t,\phi _{s}^{R}(r))+\mathsf{a}(\phi _{s},r),
\]
where $t,s\in \Bbb{T}$, $r\in R$.

In what follows we write simply $\phi _{t}$ instead of
$\phi_{t}^{R}$.

\subsection{$\frak{a}$-cocycle over flag bundles}

The cocycle $\mathsf{a}$ over $R$ can be factored to a cocycle
over the flag bundle $\Bbb{F}Q$, which is the associated bundle
$Q\times _{G}{\Bbb F}$ obtained by the left action of $G$ on
${\Bbb F}$. The construction is as follows:

Take the closed subgroup $MN$ of $G$. The quotient space $Q/MN$
can be identified with the associated bundle $Q\times _{G}G/MN$
(see \cite {kn}, Proposition 5.5). Let $P=MAN$ be the minimal
parabolic subgroup, so that the flag bundle $\Bbb{F}Q\simeq Q/P$.
Then there is a natural fibration $Q/MN\rightarrow Q/P=\Bbb{F}Q$,
$q\cdot MN\mapsto q\cdot P$, $q\in Q$, whose fiber is $P/MN\simeq
A$. Since $MN$ is normal in $P$, it follows that this fibration is
a principal bundle over $\Bbb{F}Q$ with structural group  $A $.
This bundle is trivializable (because $A$ is diffeomorphic to an
Euclidian space, see \cite{kn}). An explicit global cross section
is given by a $K$-reduction $R$ of $Q$. In fact, any element of
$\Bbb{F}Q$ can be written $r\cdot b_{0}$, $r\in R$, where $b_{0}$
is the origin of $\Bbb{F}$. Consider the map
\[
\chi (r\cdot b_{0})=r\cdot MN\in Q/MN.
\]
It is well defined because $r_{1}\cdot b_{0}=r_{2}\cdot b_{0}$
implies that  $r_{2}=r_{2}\cdot m$, $m\in M$, and hence
$r_{1}\cdot MN=r_{2}\cdot MN$. It is clearly a cross section of
$Q/MN\rightarrow \Bbb{F}Q$. Its continuity follows from the equality
 $\chi(q \cdot b_0 ) = {\sf R}(q) \cdot MN$.

Now if $\phi _{t}$ is a flow of automorphisms of $Q$, $t\in
\Bbb{T}$, then
the cross section $\chi $ defines a continuous $\frak{a}$-valued cocycle $\mathsf{a}:%
\Bbb{T}\times \Bbb{F}Q\to \frak{a}$ by $\mathsf{a}(t,\xi )=\log a_{t}$ where

\[
\phi _{t}(\chi (\xi ))=\chi (\phi _{t}(\xi ) )\cdot a_{t},\quad \xi
\in \Bbb{F}Q
\]
(cf. Section \ref{secmorgeral}). This is essentially the cocycle
over $R$ defined from the Iwasawa decomposition (which justifies
our use of the same notation). In fact, let $\xi =r\cdot b_{0}\in
\Bbb{F}Q$ where $r \in R$ and $b_0$ is the origin of $\Bbb{F}$.
Take the Iwasawa decomposition $\phi _{t}(r)=r_{t}\cdot
a_{t}n_{t}\in R\cdot AN$. Then $\phi _{t}(\xi )=r_{t}b_{0}$ so
that $\chi (\phi _{t}(\xi ) )=r_{t}\cdot MN$.  Also, since
$n_{t}\in N$ fixes $MN$, we have $\phi _{t}(\chi (\xi
))=r_{t}\cdot a_{t}MN = \chi(\phi _{t}(\xi )) \cdot a_t$, in terms
of the right action of $A$ on $Q/MN$. This means that
\[
{\sf a}(t, r\cdot b_0) = \log a_t = {\sf a}(t,r).
\]

\subsection{Cartan and polar decompositions\label{sec:cartan-polar}}

Fix a Cartan decomposition $G=KS$ with $K$ as in the Iwasawa decomposition
and the corresponding decomposition $\frak{g}=\frak{k}\oplus \frak{s}$ of
the Lie algebra $\frak{g}$ of $G$. As before let $R$ be a $K$-reduction of $%
Q $. Then for every $q\in Q$ there exists $g\in G$ and $r^{\prime
}\in R$ such that $q=r^{\prime }\cdot g$.  The $K$-component of the
Cartan decomposition $g=ks$ of $g$ leaves $R$ invariant so that, by
the uniqueness of the Cartan decomposition in $G$, it follows that
every $q \in Q$ can be written uniquely as
\[
q=r\cdot s,\quad r\in R,\quad s\in S
\]
which exhibits the bundle $Q$ as the product $Q \approx  R\times S$. We let
\[
\mathsf{R}^{C}:Q\to R,\quad q\mapsto r,\qquad \mathsf{S}:Q\to S,\quad
q\mapsto s
\]
be the associated projections. (We emphasize that the two projections $%
\mathsf{R}$ and $\mathsf{R}^{C}$ against the Iwasawa and Cartan
decompositions, respectively, are not equal.)  By the continuity of
the Cartan decomposition of $G$ and the local triviality of $Q$, it
follows that these projections are continuous. Moreover they satisfy
the following properties:

\begin{enumerate}
\item  $\mathsf{R}^{C}(q\cdot k)=\mathsf{R}^{C}(q)\cdot k$ and $\mathsf{S}%
(q\cdot k)=k^{-1}\mathsf{S}(q)k$ if $q\in Q$ and $k\in K$. (This is a
consequence of the fact that $kSk^{-1}=S$.)

\item  $\mathsf{R}^{C}(r)=r$ and $\mathsf{S}(r\cdot g)=\mathsf{S}(g)$ if $%
r\in R$, $g\in G$, where we denote also by $\mathsf{S}(g)$ the $S$-component
of $g\in G=K\times S$.

\item  $\mathsf{S}(q\cdot g)=\mathsf{S}(\mathsf{S}(q)g)$ if $q\in Q$ and $%
g\in G$. (In fact, $q\cdot g=\mathsf{R}(q)\cdot \mathsf{S}(q)g$ so that if $%
\mathsf{S}(q)g=kt$ is the Cartan decomposition of $\mathsf{S}(q)g$ in $G$
where $t=\mathsf{S}(\mathsf{S}(q)g)$ then $q\cdot g=\mathsf{R}(q)\cdot kt$.
Hence $\mathsf{S}(q\cdot g)=t=\mathsf{S}(\mathsf{S}(q)g)$.)
\end{enumerate}

Now fix a Weyl chamber $A^{+}$ sitting inside a maximal abelian
$A\subset S$ and consider the polar decomposition
$G=K(\mathrm{cl}A^{+})K$. If $g=ks\in KS $ then $g=uhv$ where
$u=kv^{-1}\in K$ and $s=v^{-1}hv$, $h\in \mathrm{cl}A^{+}$.
Combining this polar decomposition with the Cartan decomposition
of $Q$ we can write every $q\in Q$ as
\[
q=r\cdot hv,\quad r\in R,\quad h\in \mathrm{cl}A^{+},\quad v\in K.
\]
Here the component $h\in \mathrm{cl}A^{+}$ is uniquely defined (although $r$
and $v$ are not). Thus we have a well defined map
\[
\mathsf{A}^{+}:Q\to \mathrm{cl}A^{+},\quad q\mapsto h
\]
which satisfies $\mathsf{A}^{+}(q\cdot k)=\mathsf{A}^{+}(q)$ if $q\in Q$, $%
k\in K$.

In the sequel we denote with the corresponding lower case letters the
logarithms of the above maps:
\[
\mathsf{s}(q)=\log \mathsf{S}(q)\in \frak{s}\qquad \mathsf{a}^{+}(q)=\log
\mathsf{A}^{+}(q)\in \mathrm{cl}{\frak{a}}^{+}.
\]

If one starts instead with a $Z_{\Theta}$-principal bundle
$Q_{\Theta}$, using the Cartan and Iwasawa decompositions of
$Z_{\Theta}$ and the same arguments as above one gets the
existence of a $K_{\Theta}$-reduction $R_{\Theta}$ of
$Q_{\Theta}$, a Cartan decomposition of $q \in Q_{\Theta}$ given
by
\[
  q = rs, \quad r \in R_{\Theta}, \quad s \in S\cap Z_{\Theta},
\]
and an Iwasawa decomposition of $q \in Q_{\Theta}$ given by
\[
q = rhn, \quad r \in R_{\Theta}, \quad  h \in A,\, n \in
N({\Theta}).
\]
Also, if $Q_{\Theta}$ is a $Z_{\Theta}$-reduction of the
$G$-principal $Q$, then one can choose a $K$-reduction $R$ of $Q$
so that $R_{\Theta} \subset R$.

\subsection{Reductive groups\label{sec-reductive}}

The above cocycles can be also be defined for flows evolving on
principal bundles with reductive Lie groups. This is because the Lie algebra $\overline{\frak{g%
}}$ of a reductive Lie group $\overline{G}$ decomposes as $\overline{\frak{g}%
}=$ $\frak{g}\oplus \frak{z}$, with $\frak{g}$ semi-simple and
$\frak{z}$ the center. Hence one can add a central component to
$\frak{a} $ (that is, the split part of $\frak{z}$) and get a
cocycle on a larger vector space. Also, under natural conditions
there is available Cartan and Iwasawa decompositions of
$\overline{G}$ as well as parabolic subgroups and flag manifolds,
allowing analogous decompositions of $\overline{G}$-principal
bundles. We refer to Knapp \cite{kapp}, Chapter VII, for the theory
of reductive Lie groups.

We will not pursue here the detailed extension of the
decompositions for reductive principal bundles. Instead we draw
some comments regarding the group $\mathrm{Gl}\left(
n,\Bbb{R}\right) $, which show up in the vector bundles (see
Section \ref{secvecbundl}). This group is reductive with Iwasawa
decomposition \textrm{O}$\left( n,\Bbb{R}\right) \overline{A}N$,
where $N$ is the unipotent upper triangular group and
$\overline{A}$ the group of diagonal matrices with positive
entries. The Lie algebra $\overline{\frak{a}}$ of $\overline{A}$
decomposes as $\overline{\frak{a}}=\Bbb{R}\cdot \mathrm{id}\oplus
\frak{a}$ where $\frak{a}$ is the subspace of zero trace diagonal
matrices. Thus the $\overline{\frak{a}}$-valued cocycle $\overline{\mathsf{a}%
}(t,\xi)$ (constructed the same way as above) writes $\overline{\mathsf{a}}%
(t,\xi)=\mathsf{a}_{0}(t,\xi )\mathrm{id}+\mathsf{a}(t,\xi )$,
$\xi \in \Bbb{F}Q$, where $\mathsf{a}_{0}(t,\xi )$ is a real
valued cocycle which is essentially a trace of a matrix. Actually
$\mathsf{a}_{0}$ is constant as a function  of $\xi $ because
scalar matrices belong to the center of $G$.  Finally the action
of $\mathrm{Gl}\left( n,\Bbb{R}\right) $ factors to an action of
$\mathrm{Sl}\left( n,\Bbb{R}\right) $ (see \cite{smbflow}), so the
dynamics on the flag bundle is the same as for a flow on a
$\mathrm{Sl}(n,\Bbb{R})$-bundle.

\subsection{Spectra}

In the sequel we denote by $\Lambda _{\mathrm{Mo}}\left(
\mathcal{M}\right) $ and $\Lambda _{\mathrm{Ly}}\left(
\mathcal{M}\right) $ the Morse and Lyapunov vector spectra for the
$\frak{a}$-cocycle over the flag bundle $\Bbb{F}Q$ defined via the
trivializable principal bundle $Q/MN\rightarrow \Bbb{F}Q$. In
these expressions $\mathcal{M}$ is a chain component of the flow
on $\Bbb{F}Q$.

For latter reference we note that these spectra does not depend on
the trivialization of $Q/MN\rightarrow \Bbb{F}Q$ (and of the
$K$-reduction of $Q$). This is because two different
trivializations yield cohomologous cocycles, which in turn have
the same spectra, since the base is compact.

\section{Dynamics on flag bundles and block reduction\label{secpartydyn}}

In this section we describe a $\phi _{t}$-invariant subbundle
$Q_{\phi }\subset Q\rightarrow X$ whose structure group is the
centralizer $Z_{H_{\phi }}$ of an element $H_{\phi }\in \frak{a}$
(or of its exponential $h_{\phi }=\exp H_{\phi }\in A$). In case
$G$ is a linear group $H_{\phi }$ is a diagonal matrix and the
centralizer $Z_{H_{\phi }}$ is a subgroup of block diagonal
matrices (if the real eigenvalues of $H_{\phi }$ are ordered
decreasingly). By analogy to this case we call $Q_{\phi }$ the
block reduction and $H_{\phi }$ the block form of $\phi _{t}$.

The invariant subbundle $Q_{\phi }$ is constructed after the
results of \cite {smbflow}, \cite{msm} that give the chain
recurrent components of $\phi _{t}$ on the flag bundle $\Bbb{F}Q$.
These chain components in turn are modelled fiberwise by the
dynamics of the flow on the flag manifold ${\mathbb F}$ defined by
a one-parameter subgroup $\exp tY$ of $G$. The description of this
dynamics can be found \cite{dkv}, Section 3. In the next
paragraphs we recall its main features and establish related
notation. For the details and proofs see \cite{dkv}.

Let $\frak{g}=\frak{k}\oplus \frak{s}$ be a Cartan decomposition and $\frak{a%
}\subset \frak{s}$ a maximal abelian subspace. Then for every $H\in
\frak{a}$ the vector field $\widetilde{H}$ induced on $\Bbb{F}$ has
flow $\exp tH$ and is the gradient of a function with respect to a
Riemannian metric (Borel metric) in $\Bbb{F}$. Let $Z_{H}=\{g\in
G:\mathrm{Ad}\left( g\right) H=H\}$
be centralizer of $H$ in $G$ and put $K_{H}=Z_{H}\cap K$. Let $\frak{a}%
^{+}\subset \frak{a}$ be a Weyl chamber and suppose that $H\in \mathrm{cl}%
\frak{a}^{+}$. If $b_{0}$ denotes the origin of $\Bbb{F}$
(corresponding to $\frak{a}^{+}$) then the connected components of
the singularity set of $\widetilde{H}$ are the orbits $K_{H}\cdot
wb_{0}$ with $w$ running through the Weyl group $\mathcal{W}$,
which coincide with $Z_{H} \cdot wb_{0}$. These singularities are
the fixed points  the one-parameter subgroup $\exp(tH)$, $t \in
\Bbb{R}$.  These singularities are transversally hyperbolic
submanifolds of $\Bbb{F}$ (that is, $\widetilde{H}$ has Morse-Bott
dynamics) where $K_{H}\cdot b_{0}$ is the only attractor and
$K_{H}\cdot w_{0}b_{0}$ is the only repeller, where $w_{0}$ is the
element of $\mathcal{W}$ with largest length (principal
involution). Also, $K_{H}\cdot w_{1}b_{0}=K_{H}\cdot w_{2}b_{0}$
if and only if $w_{1}w_{2}^{-1}H=H$, that is, if and only if the
cosets $\mathcal{W}_{H}w_{1}$ and $\mathcal{W}_{H}w_{2}$ are
equal, where $\mathcal{W}_{H}$ is the subgroup of $\mathcal{W}$
fixing $H$.  In particular, if $H\in \frak{a}^{+}$ is regular then
there are $\left| \mathcal{W}\right| $ isolated singularities.

The stable (respectively unstable) manifold of the component
$K_{H}\cdot wb_{0}$ is the orbit $P_{H}^{-}\cdot
wb_{0}=N_{H}^{-}K_{H}\cdot wb_{0}$ (respectively $P_{H}^{+}\cdot
wb_{0}=N_{H}^{+}K_{H}\cdot wb_{0}$) of the parabolic subgroup
$P_{H}^{-}$ (respectively $P_{H}^{+}$) which is the normalizer of
the parabolic subalgebra $\frak{p}_{H}^{-}$ defined as the sum of
the eigenspaces of $\mathrm{ad}\left( H\right) $ associated with
the eigenvalues $\leq 0$ (respectively $\geq 0$).

If $Y=\mathrm{Ad}\left( g\right) H$, $H\in
\mathrm{cl}\frak{a}^{+}$, then $g$ conjugates the induced vector
fields $\widetilde{H}$ and $\widetilde{Y}$, so that the
singularities of $\widetilde{Y}$ are $g\left( K_{H}\cdot
wb_{0}\right) $ with $g\left( K_{H}\cdot b_{0}\right) $ the
attractor and $g\left( K_{H}\cdot w_{0}b_{0}\right) $ the
repeller. Their stable and unstable manifolds are respectively
$g\left( P_{H}^{\mp }\cdot wb_{0}\right) $. In the sequel we write
$\mathrm{fix}\left( Y,w\right) $ for $g\left( K_{H}\cdot
wb_{0}\right) $, the  $w$-fixed points of $Y$, and
$\mathrm{st}\left( Y,w\right) $ for the corresponding stable
manifold $g\left( P_{H}^- \cdot wb_{0}\right)$.

\begin{exemplo}
If $G=\mathrm{Sl}\left( n,\Bbb{R}\right) $ then $\Bbb{F}$ is the
manifold of complete flags $(V_{1}\subset$$\cdots$$\subset V_{n-1})$
of subspaces $V_{i}\subset \Bbb{R}^{n}$ with $\dim V_{i}=i$. We can take $\frak{%
a}$ to be the space zero trace of diagonal matrices and
$\frak{a}^{+}$ the cone of the matrices whose eigenvalues are
ordered decreasingly. If $H\in \frak{a}^{+}$ then $\widetilde{H}$
has isolated singularities which are the flags spanned by
eigenvalues of $H$. The attractor is $b_{0}=\left( \langle
e_{1}\rangle \subset \langle e_{1},e_{2}\rangle \subset \cdots
\right) $, and if $w\in \mathcal{W}$ is a permutation then
$\mathrm{fix}\left( H,w\right) $ is obtained by $w$-permuting the
basic vectors from $b_{0}$. A nonregular $H\in
\mathrm{cl}\frak{a}^{+}$ has repeated eigenvalues which appear in
blocks, and $Z_{H}$ is the subgroup of block diagonal matrices,
where the blocks have the same size as those of $H$. $K_{H}$ is the
subgroup of orthogonal matrices in $Z_{H}$. Then it is easy to write
down the fixed
point components $K_{H}\cdot wb_{0}$. For instance the attractor $\mathrm{fix%
}\left( H,1\right) $ is made of flags adapted to the block decomposition of $%
H$.
\end{exemplo}

Let $\pi :Q\to X$ be a principal bundle with semi-simple
structural group $G$ and compact Hausdorff base space  $X$. Let
$\phi _{t}$ be a right invariant flow on $Q$ which is chain
transitive on $X$. The flow $\phi _{t}$ induces a flow in the flag
bundle $\Bbb{F}Q$. The following theorem of \cite{smbflow}, \cite
{msm} gives the chain components of this induced flow.

\begin{teorema}
\label{teo:tipo-parabol-flow}The flow $\phi _{t}$ on the maximal
flag bundle ${\Bbb{F}}Q\to X$ admits a finest Morse decomposition,
whose Morse components are determined as follows: There exists
$H_{\phi }\in \mathrm{cl}\frak{a}^{+}$ and a $\phi _{t}$-invariant
map
\[
f_{\phi }:Q\to \mathrm{Ad}(G)H_{\phi }\qquad f_{\phi }\left( \phi _{t}\left(
q\right) \right) =f_{\phi }\left( q\right)
\]
into the adjoint orbit of $H_{\phi }$, which is equivariant, that
is, $f_{\phi }\left( q\cdot g\right) =\mathrm{Ad}\left(
g^{-1}\right) f_{\phi }\left( q\right) $, $q\in Q$, $g\in G$. The
chain components $\mathcal{M}\left(
w\right) $ are parameterized by the Weyl group $\mathcal{W}$ and each $%
\mathcal{M}(w)$ is given fiberwise as the fixed point set
\[
{\mathcal{M}}(w)_{\pi (q)}=q\cdot \mathrm{fix}(f_{\phi }(q),w),\quad q\in Q.
\]
There is just one attractor component
$\mathcal{M}^{+}=\mathcal{M}(1)$ and only one repeller
$\mathcal{M}^{-}=\mathcal{M}(w_{0})$ and their dynamical ordering
is the reverse of the algebraic Bruhat-Chevalley order of
$\mathcal{W}$.
\end{teorema}

In this theorem the case $H_{\phi }=0$ is not ruled out. In this case the
flow is chain transitive on the flag bundles.

In a partial flag bundle ${\Bbb{F}}_{\Theta }Q$, $\Theta \subset \Sigma $,
there exists also the finest Morse decomposition, whose components $\mathcal{%
M}_{\Theta }\left( w\right) $ are the projections of the components $%
\mathcal{M}\left( w\right) \subset \Bbb{F}Q$. These projections are also
given fiberwise as fixed points of $f_{\phi }\left( q\right) $, $q\in Q$.

Let $\Theta _{\phi }=\{{\alpha \in \Sigma :\alpha \left( H_{\phi
}\right) =0\}}$. Then the corresponding flag bundle
$\Bbb{F}_{\Theta _{\phi }}Q$ has
special properties. Namely, the attractor chain component $\mathcal{M}%
_{\Theta _{\phi }}^{+}=\mathcal{M}_{\Theta _{\phi }}\left(
1\right) $ meets
each fiber in a singleton, because the attractor fixed point component of
$H_{\phi }$ on $\Bbb{F}_{\Theta _{\phi }}$ reduces to a point. The flag $\Bbb{%
F}_{\Theta _{\phi }}$ (or the corresponding parabolic subgroup
$P_{\Theta
_{\phi }}$, or else $\Theta _{\phi }$) were called the parabolic type of $%
\phi $ in \cite{smbflow}, \cite{msm}. The parabolic type of the
reversed flow is the dual $\Theta _{\phi }^{*}$ (see
\cite{smbflow}) and in the flag
bundle $\Bbb{F}_{\Theta _{\phi }^{*}}Q$ the repeller component $\mathcal{M}%
_{\Theta _{\phi }^{*}}^{-}$ meets the fibers in singletons.
Furthermore the centralizer $Z_{H_{\phi }}=P_{\Theta _{\phi }}\cap
P_{\Theta _{\phi }}^{-}$. In the sequel we say that $H_{\phi }$ is
the block form of $\phi$.

Now let $Q_{\phi }=f_{\phi }^{-1}\left( H_{\phi }\right) \subset Q$. Since
the adjoint orbit \textrm{Ad}$\left( G\right) H_{\phi }$ is identified to
the homogeneous space $G/Z_{H_{\phi }}$, it follow by a well known fact that
$Q_{\phi }$ is a subbundle of $Q$ with structural group $Z_{H_{\phi }}$ (see
Kobayashi-Nomizu \cite{kn}, Proposition I.5.6, where the equivariant map $%
f_{\phi }$ must be viewed as a section of the associated bundle $Q\times
_{G}G/Z_{H_{\phi }}$). The $\phi _{t}$-invariance of $f_{\phi }$ implies the
invariance of $Q_{\phi }$ under the flow.

\begin{definicao}
If $H_{\phi }$ is the block form of the flow then the $Z_{H_{\phi }}$%
-subbundle $Q_{\phi }$ is called the block reduction of $\phi $.
\end{definicao}

\begin{remar}
If $Q\simeq X\times G$ is a trivial bundle then the flow is
defined by a cocycle $\rho \left( t,x\right) $ with values in $G$:
$\phi _{t}\left( x,g\right) =\left( t\cdot x,\rho \left(
t,x\right) g\right) $. By choosing another trivialization of $Q$
the cocycle is changed by cohomologous one. Now the block
reduction $Q_{\phi }$ is a subbundle of a trivial bundle. It may
not be trivial. But if this happens then the original cocycle is
cohomologous to a cocycle taking values in the subgroup
$Z_{H_{\phi }}$.
\end{remar}

\begin{remar}
The subgroup $Z_{H_{\phi }}$ is algebraic, hence it contains the
algebraic hull of the flow (see Feres \cite{fer}, Zimmer \cite{z}).
\end{remar}

The centralizer $Z_{H_{\phi }}$ is a reductive Lie subgroup of
$G$. As such it has an Iwasawa $Z_{H_{\phi }}=K_{H_{\phi }} A
N_{H_{\phi }}$, with compact component $K_{H_{\phi }}$ (see Knapp
\cite{kapp}, Chapter VII). Hence, as observed in Section
\ref{sec:cartan-polar}, the $Z_{H_{\phi }}$-bundle admits a
reduction to the compact group $K_{H_{\phi }}$.  We denote this
subbundle by $R_{\phi }$.

The following statement expresses the chain components on the flag
bundles in terms of the block reduction.

\begin{proposicao}
\label{propchaincompdynreduct} For $\Theta \subset \Sigma$, let
$\mathcal{M}_{\Theta }(w)$, $w\in \mathcal{W}$, be a chain
component of $\phi _{t}$ on the flag bundle $\Bbb{F}_{\Theta }Q$.
Then
\[
{\mathcal{M}}_{\Theta }(w)=\{q\cdot wb_{\Theta }:\,q\in Q_{\phi }\}=\{r\cdot
wb_{\Theta }:\,r\in R_{\phi }\},
\]
where $b_{\Theta }$ is the origin of $\Bbb{F}_{\Theta }$.
\end{proposicao}

\begin{demonstracao}
In fact, if $q\in Q_{\phi }$ then ${\mathcal{M}}_{\Theta }(w)_{\pi
(q)}=q\cdot \mathrm{fix}(H_{\phi },w)$, but $\mathrm{fix}(H_{\phi },w)=Z_{H_{\phi }}wb_{\Theta }$,
 because $H_{\phi }\in \mathrm{cl}\frak{a}^{+}$. Then the first equality follows since
 $Z_{H_{\phi }}$ acts transitively on the right in the fiber of $q$ in $Q_{\phi }$.
  The proof of the second equality is similar, by taking $K_{H_{\phi }}$ and $R_{\phi }$
instead of $Z_{H_{\phi }}$ and $Q_{\phi }$.
\end{demonstracao}

We conclude this section with the construction of a linearization
around the attractor of the induced flow in the flag bundle
$\Bbb{F}_{\Theta ({\phi}) }Q$. This linearization is a flow which
evolves on a subbundle of the  tangent bundle to the fibers of
$\Bbb{F}_{\Theta }Q$. It will be the main technical tool in the
proof of the estimates of Section \ref {seclocal}. More details
about this construction can be found in the forthcoming paper
\cite{psmsec}.

We write $\Theta = \Theta \left( \phi \right) =\{\alpha \in \Sigma
:\alpha \left( H_{\phi }\right) =0\}$ for the parabolic type of
$\phi $. As before $Z_{\Theta }$ stands for the centralizer of
$H_\phi$ and $Q_{\phi}$ the $Z_{\Theta }$-block reduction. Let
$\frak{n}_{\Theta }^{-}$ be the nilradical of $\frak{p}_{\Theta
}^{-}$ (which identifies to the tangent space at the origin of
$\Bbb{F}_{\Theta }$). The group $Z_{\Theta }$ normalizes
$\frak{n}_{\Theta }^{-}$, and hence acts linearly on
$\frak{n}_{\Theta }^{-}$ by the adjoint representation. Hence, we
can build the associated bundle
\[
\mathcal{V}_{\Theta }=Q_\phi\times _{Z_{\Theta }}\frak{n}_{\Theta
}^{-}\rightarrow X.
\]

Since the $Z_{\Theta }$-action is linear, it follows that
$\mathcal{V}_{\Theta }\rightarrow X$ is a vector bundle and the
flow $\Phi _{t}$ induced by $\phi _{t}$ on $\mathcal{V}_{\Theta }$
is linear.

Let $b_{\Theta }$ be the origin in $\Bbb{F}_{\Theta }$ and define the subset
\[
\Bbb{B}_{\Theta }=Q_\phi\cdot N_{\Theta }^{-}b_{\Theta }
\]
and the mapping $\Psi :\mathcal{V}_{\Theta }\rightarrow
\Bbb{B}_{\Theta }$,
\[
\Psi \left( q\cdot X\right) =q\cdot \left( \exp X\right) b_{\Theta
}\qquad q\in Q_\phi,X\in \frak{n}_{\Theta }^{-}.
\]
Note that $\Psi $ is well defined because if $q^{\prime }\cdot
X^{\prime }=q\cdot X$ then there exists $g\in Z_{\Theta }$ with
$q^{\prime }=q\cdot g$ so that $X^{\prime }=\mathrm{Ad}\left(
g\right) X$. But $\exp X^{\prime }=g\exp Xg^{-1}$ and $gb_{\Theta
}=b_{\Theta }$ so that $q^{\prime }\cdot \left( \exp X^{\prime
}\right) b_{\Theta }$ coincides with $q\cdot \left( \exp X\right)
b_{\Theta }$. Moreover, $\Psi $ is a homeomorphism because the
mapping $X\in \frak{n}_{\Theta }^{-}\mapsto \left( \exp X\right)
b_{\Theta }\in N_{\Theta }^{-}b_{\Theta }$ is a homeomorphism (by
Proposition 3.6  of \cite{dkv}).

\begin{proposicao}
\label{proplinvecvteta}The following statements are true.

\begin{enumerate}
\item  $\Bbb{B}_{\Theta }$ is an open and dense $\phi _{t}$-invariant subset
of $\Bbb{E}_{\Theta }$ which contains the attractor component $\mathcal{M}%
_{\Theta }^{+}=\Psi \left( \mathcal{V}_{\Theta }^{0}\right) $,
where $\mathcal{V}_{\Theta }^{0}$ is the zero section of
$\mathcal{V}_{\Theta }$.

\item  $\phi _{t}$ and $\Phi _{t}$ are conjugate under $\Psi $:
\[
\phi _{t}\left( \Psi \left( v\right) \right) =\Psi \left( \Phi
_{t}\left( v\right) \right), \qquad v\in \mathcal{V}_{\Theta }.
\]
\end{enumerate}
\end{proposicao}

\begin{demonstracao}
The $\phi _{t}$-invariance of $\Bbb{B}_{\Theta }$ follows from the
$\phi _{t} $-invariance of $Q_{\phi}$. To see that it is open and
dense we note that $N_{\Theta }^{-}b_{\Theta }$ is the open Bruhat
cell in $\Bbb{F}_{\Theta }$, hence $q\cdot N_{\Theta }^{-}b_{\Theta
}$ is open and dense in its fiber. Making $q$ run through
$Q_{\phi}$, the result follows. Now,
\[
\mathcal{M}_{\Theta }^{+}=Q_\phi\cdot b_{\Theta }=\Psi \left( \mathcal{V%
}_{\Theta }^{0}\right) \subset \Bbb{B}_{\Theta },
\]
concluding the proof of the first statement.  To see the conjugacy
take $v=q\cdot X$, $q\in Q_\phi$, $X\in \frak{n}_{\Theta }^{-}$.
Then by definition $\Phi _{t}\left( v\right) =\phi _{t}\left(
q\right) \cdot v$, so that
\[
\Psi \left( \Phi _{t}\left( v\right) \right) =\phi _{t}\left(
q\right) \cdot \left( \exp X\right) b_{\Theta }=\phi _{t}\left(
q\cdot \exp X\right) b_{\Theta }=\phi _{t}\left( \Psi \left(
v\right) \right) ,
\]
concluding the proof.
\end{demonstracao}

Finally we endow $\mathcal{V}_{\Theta }\rightarrow X$ with a
natural metric $\left( \cdot ,\cdot \right) $ given by
\[
\left( r\cdot X,\,r\cdot Y\right) =B_{\theta }(X,Y),\quad r\in R_{\Theta
},\,X,Y\in \frak{n}_{\Theta }^{-}
\]
where $B_{\theta }\left( \cdot ,\cdot \right) $ is the inner
product in the Lie algebra defined by the Cartan involution
$\theta $. That this in fact defines a metric in the whole
$\mathcal{V}_{\Theta }$ follows from the Iwasawa decomposition
$Q_\phi = R_{\phi} A N(\Theta)$, where $A N(\Theta)$ normalizes
$\frak{n}_\Theta^-$.

\section{$\mathcal{W}$-invariance of Lyapunov exponents}

The purpose of this section is to prove the invariance of the set
of $\frak{a}$-Lyapunov exponents under the Weyl group (see
Corollary \ref{corinvariance}). The proof is based on the concept
of regular sequences in a symmetric space (cf.  \cite{ka}).

\subsection{Regular sequences in $G$}

We discuss here some asymptotic properties of sequences in $G$ that will be
used afterwards in the proof of the invariance of the $\frak{a}$-Lyapunov
exponents under the Weyl group. Most of the results here are taken from \cite
{ka}, although we need to adapt them to our $\frak{a}$-valued exponents.

For a sequence $g_{k}\in G$, $k\in \Bbb{Z}^{+}$, there are two
related limits in $\frak{a}$ that show up in the Iwasawa or Cartan
decompositions of $g_{k}$.

\begin{enumerate}
\item  The Lyapunov exponent of $g_k$ at $b\in \Bbb{F}$:
\[
\lambda \left( g_{k},b\right) =\lim_{k\to +\infty }\frac{1}{k}\log
{\sf a}(g_ku),
\]
where $b = ub_0$, $u \in G$. This limit depends only on $b$
because if $b = u'b_0$, $u' \in G$, then $u' = up$, for $p \in P$,
so by Equation (\ref{eq-propr-aditiva-iwasawa}) one has ${\sf a}(g_ku') =
{\sf a}(g_ku) + {\sf a}(p)$, and in the limit the term ${\sf
a}(p)$ disappears.

\item  The polar exponent of $g_k$:
\[
\lambda ^{+}\left( g_{k}\right) =\lim_{k\to +\infty }\frac{1}{k}\log
{\sf a}^+(g_k),
\]
which depends on the sequence only, contrary to the Lyapunov
exponents.
\end{enumerate}

Of course  these limits may not exist.  Next consider the left
coset symmetric space $K\backslash G$ (we use this form instead of
$G/K$ to match to the Iwasawa decomposition $G=KAN$ as well as to
the right action of $G$ on the principal bundle $Q\rightarrow X$).
Let $d$ be the $G$-invariant distance in $K\backslash G$, which is
uniquely determined by
\[
d(x_{0}\cdot \exp (X),x_{0})=|X|_{\theta },\quad X\in \frak{s}
\]
where $x_{0}$ is the origin of $K\backslash G$.

Following \cite{ka} a sequence $g_{k}\in G$ is said to be regular if
there exists $D\in \frak{s}$ such that $d\left( x_{0}\cdot
g_{k},x_{0}\cdot \exp kD\right) $ has sublinear growth as
$k\rightarrow +\infty $, that is, if
\begin{equation}
\frac{1}{k}d\left( x_{0}\cdot g_{k},x_{0}\cdot \exp kD\right) \rightarrow 0.
\label{fordefregular}
\end{equation}
In this case we say that $D$ is the asymptotic ray of $g_{k}$. (More
precisely, $x_0 \cdot g_{k}$ is asymptotic to the geodesic ray $x_0
\cdot \exp tD$, $t>0$.)  If $u\in K$ it follows that the sequence $g_{k}u$ is
also regular with asymptotic ray $\mathrm{Ad}\left( u^{-1}\right) D$.

Write the polar decomposition of $g_{k}$ as $g_{k}=u_{k}h_{k}v_{k}
\in K\left( \mathrm{cl}A^{+}\right) K$. Then one of the main
results of \cite{ka} says that the sequence is regular if and only
if it satisfies the following two conditions (see \cite{ka},
Theorem 2.1):
\begin{enumerate}
\item  The polar exponent $\lambda ^{+}\left( g_{k}\right) =\lim \frac{1}{k}%
\log h_{k}$ (or equivalently the limit $h_{k}^{1/k}$) exists, and

\item  $\frac{1}{k}d(x_{0}\cdot g_{k},\,x_{0}\cdot g_{k+1})\rightarrow 0$,
that is, $x_{0}\cdot g_{k}$ has sublinear growth in $K\backslash G$.
\end{enumerate}

Moreover, in this case the asymptotic ray of $g_k$ is given by
$D=\mathrm{Ad}(u) H^{+}$ where $u \in K$ and $H^{+}=\lambda
^{+}\left( g_{k}\right) \in \mathrm{cl}\frak{a}^{+}$ is the polar
exponent of $g_k$.

The next statement gives the Lyapunov exponents of a regular sequence.
\begin{proposicao}
\label{propkaimanovich}Suppose $g_{k}$ is a regular sequence in $G$ with
asymptotic ray $D\in \frak{s}$. Then the following statements hold.

\begin{enumerate}
\item  Suppose  $D\in \mathrm{cl}\frak{a}^{+}$ and take $y\in P_{D}^{-}$.
Then the sequence $g_{k}y$ is regular with the same asymptotic ray $D$.

\item  Suppose  $D \in \frak{a}$.
Then $\lim_{k\rightarrow +\infty }\frac{1}{k}\log {\sf a}(g_k)=D$.

\item  Let $H^{+} \in \mathrm{cl} \frak{a}^+$ be the polar exponent of $g_k$.
Then the Lyapunov exponent at $b\in \mathrm{st}(D,w)$ is given by
$\lambda \left( g_{k},b\right) =w^{-1}H^{+}$.

\end{enumerate}
\end{proposicao}

\begin{demonstracao}
For the first statement we note that since the eigenvalues of $\mathrm{ad}%
\left( D\right) $ on $\frak{p}_{D}^{-}$ are $\leq 0$, it follows
from the Iwasawa decomposition of $y \in P^-_D$ that the sequence
$\left( \exp kD\right) y \left( \exp \left( -kD\right) \right) $ is
bounded in $P_{D}^{-}$ so that
\[
d(x_{0}\cdot \left( \exp kD\right) y,\,x_{0}\cdot \exp kD)=d(x_{0}\cdot
\left( \exp kD\right) y\left( \exp \left( -kD\right) \right) ,\,x_{0})
\]
is bounded as a function of $k$ (cf. \cite{gjt}, Proposition 3.9).
But
\[
d(x_{0}\cdot g_{k}y,\,x_{0}\cdot \exp kD)\leq d(x_{0}\cdot
g_{k},\,x_{0}\cdot \exp kD)+d(x_{0}\cdot \exp kD,\,x_{0}\cdot \left( \exp
kD\right) y^{-1}).
\]
So that $\frac{1}{k}d(x_{0}\cdot g_{k}y,\,x_{0}\cdot \exp kD)\rightarrow 0$,
because $g_{k}$ is asymptotic to $\exp kD$ and the last term is bounded.

For the second statement write $g_{k}=u_{k}a_{k}n_{k}\in KAN$. We
have
\[ \left| \log a_{k}-kD\right| _{\theta }=d\left(
x_{0}\cdot a_{k},x_{0}\cdot \exp kD\right)
\]
and by a well known inequality of horospherical coordinates (cf.
Corollary 1.11 of \cite{Gabriele}) the right hand side is bounded
above by
\[
d\left( x_{0}\cdot a_{k}n_{k},x_{0}\cdot \exp kD\right) =d\left(
x_{0}\cdot g_{k},x_{0}\cdot \exp kD\right).
\]
Therefore, $\frac{1}{k}\log a_{k}\rightarrow D$.

Now, recall that $\mathrm{st}(D,w)=uP_{H^{+}}^{-}wb_{0}$ where $u\in K$
satisfies $D=\mathrm{Ad}(u)H^{+}$ and $b_{0}\in \Bbb{F}$ is the origin.
Hence, for $b\in \mathrm{st}(D,w)$ there exists $y\in P_{H^{+}}^{-}$ such
that $b=uywb_{0}$. In this case $\lambda \left( g_{k},b\right)
=\lim_{k\rightarrow +\infty }\frac{1}{k}\log a_{k}$ where $%
g_{k}uyw=u_{k}^{\prime }a_{k}n_{k}\in KAN$. Since $u\in K$ the
sequence $g_{k}u$ is regular with asymptotic ray $H^{+}$, so by
the first statement the same holds for $g_{k}uy$. But $w\in K$
hence $g_{k}uyw$ is asymptotic to the ray $\mathrm{Ad}\left(
w^{-1}\right) H^{+}\in \frak{a}$. Therefore, $\lambda \left(
g_{k},b\right) =\mathrm{Ad}\left( w^{-1}\right) H^{+}$, as follows
from the second statement, concluding the proof.
\end{demonstracao}

In \cite{ka} it is also given the following characterization of
regularity in terms of the Lyapunov exponents $\lambda \left(
g_{k},b\right) $.

\begin{proposicao}
\label{proplypumtodos}The sequence $g_{k}$ is regular if and only
if it has sublinear growth and the Lyapunov exponent $\lambda
\left( g_{k},b\right) $ exists at some (and hence at all) $b\in
\Bbb{F}$.
\end{proposicao}

\begin{demonstracao}
See \cite{ka}, Theorem 2.5 and its corollary where it is proved that
$g_{k}$ is regular it it has sublinear growth and $\lim
\frac{1}{k}\log a_{k}$ exists, where $g_{k}=u_{k}^{\prime
}a_{k}n_{k}\in KAN$ is the Iwasawa decomposition. But as mentioned
above $g_{k}$ is regular if and only if $g_{k}u$, $u\in K$,
is regular. Hence the existence of $\lambda \left( g_{k},b\right) $ at some $%
b\in \Bbb{F}$ together with sublinear growth implies regularity.
\end{demonstracao}

Now we specialize the above results for sequences $g_{k}$ taking
values in a centralizer subgroup $Z_{Y}$, $Y\in \frak{a}$. By the
block reduction dicussed in Section \ref{secpartydyn} the sequences
in $G$ giving rise to our exponents will ultimately belong to a
$Z_{Y}$.

\begin{proposicao}
\label{prop:kaimanovich-em-Z}Take $Y\in \mathrm{cl}\frak{a}^{+}$
and put $\Theta =\{\alpha \in \Sigma :\alpha \left( Y\right)
=0\}$. Let $g_{k}\in Z_{Y}$ be a regular sequence in $G$. Then the
asymptotic ray of $g_{k}$ has the form $D=\mathrm{Ad}(u)H$ with
$u\in K_{\Theta }$ and $H\in \frak{a}$ such that $\alpha (H)\geq
0$ for all $\alpha \in \Theta $.
\end{proposicao}

\begin{demonstracao}
Since $Z_Y = Z_\Theta$ we may take the decomposition
$g_{k}=m_{k}g_{k}^{\prime }h_{k}$, with $m_{k}\in M$,
$g_{k}^{\prime }\in G(\Theta )$, $h_{k}=\exp H_{k}\in A_{\Theta }$, where $%
H_{k}\in \frak{a}_{\Theta }$. Let $D$ be the asymptotic ray of $g_{k}$ so
that
\[
d(x_{0}\cdot g_{k},\,x_{0}\cdot \exp kD)=d(x_{0}\cdot g_{k}^{\prime
}h_{k},\,x_{0}\cdot \exp kD).
\]
This implies that $g_{k}^{\prime }h_{k}$ is also regular in $G$. By Theorem
2.5 of \cite{ka} it follows that $g_{k}^{\prime }$ is regular in $G(\Theta )$
and the limit $\frac{1}{k}H_{k}\to \widehat{H}$ exists. Hence, there exists $%
H^{\prime }\in \mathrm{cl}\frak{a}(\Theta )^{+}$ and $u\in K(\Theta )$ such
that
\[
\frac{1}{k}d^{\prime }(x_{0}^{\prime }\cdot g_{k}^{\prime },\,x_{0}^{\prime
}\cdot \left( \exp kH^{\prime }\right) u^{-1})\rightarrow 0,
\]
where $d^{\prime }$ is the distance in the symmetric space $K(\Theta
)\backslash G(\Theta )$ with origin $x_{0}^{\prime }=K(\Theta )$.
Since the embeeding $K(\Theta )\backslash G(\Theta )\hookrightarrow
K\backslash G$ is isometric, this limit still holds with the
distance $d$ of $K\backslash G$.
Now $x_{0}\cdot g_{k}^{\prime }h_{k}=x_{0}\cdot g_{k}$, and since $%
h_{k}=\exp (H_{k})\in A_{\Theta }$ comutes with $u\in K_{\Theta }$
and $A$ is abelian we multiply by $h_k$ on both arguments of the
distance function in the previous limit to conclude that
\[
\frac{1}{k}d(x_{0}\cdot g_{k},\,x_{0}\cdot \exp (kH^{\prime
}+H_{k})u^{-1})\rightarrow 0.
\]
Put $H=H^{\prime }+\widehat{H}$. Then the above limits show that
\[
\frac{1}{k}d(x_{0}\cdot g_{k},\,x_{0}\cdot \exp k(H^{\prime }+\widehat{H}%
)u^{-1})\rightarrow 0,
\]
so that the asymptotic ray of $g_{k}$ is $D=\mathrm{Ad}(u)H$. Here $u\in
K_{\Theta }$ and $\alpha (H)=\alpha (H^{\prime })\geq 0$ for all $\alpha \in
\Theta $ as claimed.
\end{demonstracao}

\subsection{Lyapunov exponents}

Let $Q=R\cdot AN$ be an Iwasawa decomposition of the principal bundle $Q$,
where $R$ is a $K$-reduction of $Q$. As in Section \ref{secdecom} the
logarithm of the $A$-component gives rise to the $\frak{a}$-cocycle $\mathsf{%
a}(t,\xi )$ over the flag bundle $\Bbb{F}Q$ and to the
$\frak{a}$-Lyapunov exponent of the flow $\phi_t$
\[
\lambda \left( \xi \right) =\lim_{t\to +\infty }\frac{1}{t}\mathsf{a}(t,\xi
)
\]
in the direction of $\xi \in \Bbb{F}Q$.

It is useful to write down  limits at the  bundle level in terms
of sequences in the  group level and hence to obtain results about
the Lyapunov exponents of the flow from the asymptotics  in $G$.
The following statement follows directly from the definitions.

\begin{proposicao}
\label{proplypseqfib}Let $r\in R$ and write $\phi _{t}(r)=r_{t}\cdot g_{t}$,
with $r_{t}\in R$, $g_{t}\in G$, $t\in \Bbb{T}$. If $b_{0}\in \Bbb{F}$ is
the origin and $b=ub_{0}$, $u\in K$, then
\[
\lambda \left( r\cdot b\right) =\lambda (r\cdot ub_{0})=\lim_{t\to +\infty }%
\frac{1}{t}\log a_{t},\quad u\in G
\]
where $g_{t}u=u_{t}^{\prime }a_{t}n_{t}\in KAN$ is the Iwasawa decomposition.
\end{proposicao}

In other words any Lyapunov exponent of the flow is a Lyapunov exponent of a
sequence in $G$.

Let $Q=R\cdot S\simeq R\times S$ be a Cartan decomposition of the
bundle (see Section \ref{secdecom}). As before we let
$\mathsf{S}:Q\to S$ stand for
the projection onto $S\simeq K\backslash G$. For any initial condition $%
q\in Q$ we define the sequence $\mathsf{S}(\phi _{k}(q))\in S$,
$k\in \Bbb{Z}^+$.

In what follows we say that $q\in Q$ is a regular point for the
flow in case $\mathsf{S}(\phi _{k}(q))$ is a regular sequence in
$S\subset G$ in the sense of last subsection. We note that $q$ is
a regular point if and only if $q\cdot g$ is regular for every
$g\in G$. That is, regularity depends only on the fiber of $q$. In
this case we say that the base point $x=\pi \left( q\right) \in X$
is regular for the flow.

We intend to show that every $\frak{a}$-Lyapunov exponent of the
flow $\phi _{t}$ is a Lyapunov exponent of some regular sequence
$\mathsf{S}(\phi _{k}(q))$. To this purpose we check first that as
a consequence of continuity combined with compactness of the base
space $X$ the sublinear growth of $\mathsf{S}(\phi _{k}(q))$
always holds.

\begin{proposicao}
For any $q\in Q$ we have $\frac{1}{k}d(x_{0}\cdot s_{k},\,x_{0}\cdot
s_{k+1})\rightarrow 0$ where $s_{k}=\mathsf{S}(\phi _{k}(q))$.
\end{proposicao}

\begin{demonstracao}
Write $q_{k}=\phi _{k}(q)\in Q$ and let $q_{k}=r_{k}\cdot s_{k}\in R\cdot S$
be its Cartan decomposition. Clearly $q_{k+1}=\phi _{1}(q_{k})$ so that
\[
q_{k+1}=r_{k+1}\cdot s_{k+1}=\phi _{1}(q_{k})=\phi _{1}(r_{k})\cdot s_{k},
\]
by right invariance. Therefore, $s_{k+1}=\mathsf{S}(q_{k+1})=\mathsf{S}(\phi
_{1}(r_{k})\cdot s_{k})=\mathsf{S}(\mathsf{S}(\phi _{1}(r_{k}))s_{k})$ (see
Section \ref{secdecom}). Hence
\begin{eqnarray*}
d(x_{0}\cdot s_{k},\,x_{0}\cdot s_{k+1})=d(x_{0}\cdot s_{k},\,x_{0}\cdot
\mathsf{S}(\phi _{1}(r_{k}))s_{k})= \\
=d(x_{0},\,x_{0}\cdot \mathsf{S}(\phi _{1}(r_{k})))=|\log \mathsf{S}(\phi
_{1}(r_{k}))|.
\end{eqnarray*}
By compactness of the base space we have that $R$ is compact.  From
the continuity of $\phi_t$ and of the Cartan decomposition it
follows then that $d(x_{0}\cdot s_{k},\,x_{0}\cdot s_{k+1})$ is
bounded, which establishes the result.
\end{demonstracao}

Now we can prove that any Lyapunov exponent of the flow on $Q$
comes from a regular sequence in $G$.

\begin{proposicao}
\label{propkaimanovichflow}Let $r\in R$. Then the following statements are
equivalent.

\begin{enumerate}
\item  $r$ is a regular point for the flow. Denote by $D\in \frak{s}$ the
asymptotic ray and by $H^+\in \mathrm{cl} \frak{a}^+$ the polar
exponent of $\mathsf{S}(\phi _{k}(r))$.

\item  The Lyapunov exponent $\lambda (r\cdot b)$ exists in one, and hence
in any direction $r\cdot b$, $b\in \Bbb{F}$, along the fiber of
$r$. In that case $\lambda (r\cdot b)=w^{-1}H^{+}$ for any  $b\in
\mathrm{st}(D,w)$.
\end{enumerate}
\end{proposicao}

\begin{demonstracao}
Take the Cartan decomposition $\phi _{k}(r)=r_{k}\cdot s_{k}\in
R\cdot S$. Then $s_{k}=\mathsf{S}(\phi _{k}(r))$ has sublinear
growth. By Proposition \ref {proplypseqfib} the Lyapunov exponent at
$r\cdot b$, $b\in \Bbb{F}$, is the Lyapunov exponent at $b$ of the
sequence $s_{k}$. Now the result follows by Proposition
\ref{proplypumtodos} which ensures that a sequence with sublinear
growth is regular if and only if some (and hence all) Lyapunov
exponent exists.
\end{demonstracao}

It follows immediately that any Lyapunov exponent of the flow is the
Lyapunov exponent of a regular sequence.

\begin{corolario}
\label{corol:espectro-pontos-reg}Take $\xi=r\cdot b\in \Bbb{F}Q$ and
assume that its Lyapunov exponent $\lambda \left( \xi \right) $
under $\phi _{t}$ exists. Then $r$ is a regular point for the flow
and $\lambda \left( \xi \right) =\lambda \left( s_{k},b\right) $
where $s_{k}=\mathsf{S}(\phi _{k}(r))$.
\end{corolario}

By Proposition \ref{propkaimanovich} the set of
$\frak{a}$-Lyapunov exponents $\lambda \left( g_{k},b\right) $,
$b\in \Bbb{F}$, of a regular sequence $g_{k} $ is invariant under
the Weyl group $\mathcal{W}$. Therefore the above corollary has as
a consequence the following symmetry of the Lyapunov spectrum.

\begin{corolario}
\label{corinvariance}Take $\xi=r\cdot b\in \Bbb{F}Q$ and assume that
its Lyapunov exponent $\lambda \left( \xi\right) $ for the flow
$\phi _{t}$ exists. Then for every $w\in \mathcal{W}$, $w\lambda
\left( \xi\right) $ is also an $\frak{a}$-Lyapunov exponent of the
flow.
\end{corolario}

In other words this last result says that the whole set of Lyapunov
exponents is invariant under the Weyl group. Since the Morse exponents are
convex combinations of Lyapunov exponents, it follows that the same result
holds for the whole Morse spectrum of the flow. This last statement will be
clarified later with the description of the Morse spectrum of each chain
component.

\begin{remar}
The Iwasawa decomposition of $Q$ gives an additive cocycle over
$\mathbb{F} Q$.  The polar decomposition of $Q$ can be shown to
give a subaditive $\frak{a}$-cocycle over $X$ in the following
way. Put ${\sf a}^+(t,x) = {\sf a}^+(\phi_t(r))$, where $r \in R$
and $\pi(r) = x$. Let $\mu \in \frak{a}^*$ be a weight of a finite
dimensional representation of $G$, then
\[
\langle {\mu, {\sf a}^+(t+s,x)} \rangle \leq
\langle {\mu, {\sf a}^+(t, \phi_s(x))} \rangle  +
\langle {\mu, {\sf a}^+(s,x)} \rangle .
\]
As usual (see \cite{ka}) an
application of the subadditive ergodic theorem to this cocycle shows
that almost all points (with respect to an invariant measure in $X$)
are regular, thus yielding the Multiplicative Ergodic Theorem for the
$\frak{a}$-Lyapunov exponents. Since these matters are exhaustively treated
in the literature we do not exploit it any further.
\end{remar}

\section{Spectrum of the attractor component and Weyl chambers\label{seclocal}}

In this section we prove one of the main results of this paper. It
locates the Morse spectrum of the attractor chain component
$\mathcal{M}^{+} \subset {\mathbb F}Q$  within an open cone in
$\frak{a}$. This cone is read off from  the parabolic type of the
flow (see Corollary \ref{corincludecone}). The proof is based on
estimates of the linear flow on the vector bundle
$\mathcal{V}_{\Theta \left( \phi \right)}\rightarrow X$ discussed
at the end of Section \ref {secpartydyn}.

We make use of the following statement about linear flows on vector bundles.

\begin{proposicao}
\label{propos:fenichel}Let $\Phi _{t}$ be a linear flow on a
vector bundle $\mathcal{V}\to X$ with norm $\left| \cdot \right| $
over a compact Hausdorff base space. Suppose that the zero section
$\mathcal{V}^{0}$ is an attractor. Then there are positive
constants $C,\mu >0$ such that
\[
\Vert \Phi _{t}\Vert _{x}\leq Ce^{-\mu t},\qquad \quad t\geq 0,
x\in X,
\]
where $\Vert \cdot \Vert $ denotes the operator norm of $|\cdot|$.
\end{proposicao}

\begin{demonstracao}
Let $U$ be an attracting neighborhood of $\mathcal{V}^{0}$ in
$\mathcal{V}$. Let $v\in \mathcal{V}$ be arbitrary. Then
${\varepsilon }v\in U$ for sufficiently small ${\varepsilon }>0$
which implies that $|\Phi _{t}({\varepsilon }v)|={\varepsilon
}|\Phi _{t}(v)|\to 0$ when $t\to +\infty $. Hence, $|\Phi
_{t}(v)|\to 0$, when $t\to +\infty $. The result then follows by
the uniformity lemma of Fenichel \cite{fer} (see also \cite{ck},
Lemma 5.2.7, and \cite{salzh}, Theorem 2.7).
\end{demonstracao}

The constant $\mu > 0$ is called a contraction exponent of the
linear flow. By compactness of the base space it is easily seen to
be independent of the  chosen metric in $\mathcal{V}$.

We have the following estimate of the values of the roots on the
$\frak{a}$-Lyapunov exponent in $\Lambda
_{\mathrm{Ly}}(\mathcal{M}^{+})$.

\begin{teorema}
\label{teoespectrlyapposit}Let $\mu >0$ be a contraction exponent of
the linear flow $\Phi_{t}$ on $\mathcal{V}_{\Theta (\phi )}$ given
in Proposition \ref{proplinvecvteta}. Then every $\frak{a}$-Lyapunov
exponent $\lambda \in \Lambda_{\mathrm{Ly}}(\mathcal{M}^{+})$
satisfies
\begin{equation}
\alpha (\lambda )\geq \mu ,\quad \alpha \in \Pi ^{+}\setminus
\langle {\Theta (\phi )}\rangle .
\label{eq:teo:localiz-espectro-atrator}
\end{equation}
\end{teorema}

\begin{demonstracao}
The main point in the proof is to relate the operator norm of
$\Phi _{t}$ with the cocycle $\mathsf{a}(t,\xi )$ over the flag
bundle. Let $\Theta =\Theta ({\phi}) $. Take $\xi \in
\mathcal{M}^{+}$ and write $\xi =r\cdot b_0$, $r \in R_\phi$. Then
for $t\in \Bbb{T}$ we have the Iwasawa decomposition $\phi
_{t}(r)=r_{t}\cdot a_{t}n_{t}\in R_{\Theta }AN^{+}(\Theta )$ and
\[
\mathsf{a}(t,\xi )=\log a_{t}.
\]
On the other hand if $v=r\cdot Y\in \mathcal{V}_{\Theta }$ with $Y\in \frak{n%
}_{\Theta }^{-}$ then $|v|=\left| Y\right| _{\theta }$. But $\Phi
_{t}(v)=\phi _{t}(r)\cdot Y$ and hence
\[
\left| \Phi _{t}(v)\right| =\left| \phi _{t}(r)\cdot Y\right|
=\left|
r_{t}\cdot \mathrm{Ad}(a_{t}n_{t})Y\right| =\left| \mathrm{Ad}%
(a_{t}n_{t})Y\right| _{\theta }.
\]
Since $N(\Theta)$ centralizes $\frak{n}^-_\Theta$ we have
$\mathrm{Ad}(n_{t})Y = Y$ so that $\left| \Phi _{t}(v)\right|
=\left| \mathrm{Ad}(a_{t})Y\right| _{\theta }$. Therefore
\[
\Vert \Phi _{t}\Vert _{x}=\Vert
\mathrm{Ad}(a_{t})|_{\frak{n}_{\Theta }^{-}}\Vert _{\theta }.
\]
Now $\mathrm{Ad}(a_{t})|_{\frak{n}_{\Theta }^{-}}$ is positive
definite so that $\Vert \mathrm{Ad}(a_{t})|_{\frak{n}_{\Theta
}^{-}}\Vert _{\theta }$ equals its greatest eigenvalue. Since the
eigenvalues are $e^{-\alpha (\log a_{t})}$, $\alpha \in \Pi ^{+}
\setminus \langle {\Theta }\rangle $, it follows that
\[
\log \Vert \Phi _{t}\Vert _{x}\geq -\alpha (\log a_{t})=-\alpha (\mathsf{a}%
(t,\xi )),\quad \alpha \in \Pi ^{+}\setminus \langle {\Theta }\rangle .
\]
By Proposition \ref{propos:fenichel} there exists $B \in
\mathbb{R}$ such that $\alpha (\mathsf{a}(t,\xi ))\geq \mu t+B $
for all roots $\alpha \in \Pi ^{+}\setminus \langle {\Theta
}\rangle $ and $t \geq 0$.

Finally, if $\lambda \in \Lambda _{\mathrm{Ly}}(\mathcal{M}^{+})$
then $\lambda =\lambda (\xi )=\lim_{t\to \infty
}\frac{1}{t}\mathsf{a}(t,\xi )$ for some $\xi \in
\mathcal{M}^{+}$, so that $\alpha \left( \lambda \right) \geq \mu
$, concluding the proof.
\end{demonstracao}

\begin{corolario}
\label{corparabtyplypexp}Take $\lambda \in \Lambda _{\mathrm{Ly}}\left(
\mathcal{M}^{+}\right) $ and write $\Theta (\lambda )=\{\alpha \in \Pi
^{+}:\alpha \left( \lambda \right) =0\}$. Then $\Theta (\lambda )$ is
contained in the parabolic type $\Theta (\phi )$.
\end{corolario}

Since the Morse spectrum $\Lambda _{\mathrm{Mo}}\left( \mathcal{M}%
^{+}\right) $ is the convex closure of $\Lambda _{\mathrm{Ly}}\left(
\mathcal{M}^{+}\right) $ it follows at once that the same estimate of the
above theorem holds for the Morse exponents.

\begin{corolario}
\label{corraizpositilambdamo}If $\lambda \in \Lambda _{\mathrm{Mo}}\left(
\mathcal{M}^{+}\right) $ then $\alpha (\lambda )>\mu >0$ for every $\alpha
\in \Pi ^{+}\setminus \langle {\Theta (\phi )}\rangle $.
\end{corolario}

Now we have the following lemma on root systems, which might be well known.
By the lack of a reference we present a proof of it.

\begin{lema}
The set $\{H\in \frak{a}:\alpha \left( H\right) >0,\alpha \in \Pi
^{+}\setminus \langle {\Theta }\rangle \}$ is the open convex cone $\mathrm{%
int}\left( \mathcal{W}_{\Theta }\mathrm{cl}\frak{a}^{+}\right) $.
Also, two cones $\mathrm{int}\left( w_i \mathcal{W}_{\Theta
}\mathrm{cl}\frak{a}^{+}\right)$, $w_1, w_2 \in \mathcal{W}$, are
either equal or disjoint.
\end{lema}

\begin{demonstracao}
Let $\mathcal{C}_{\Theta }$ be the cone of those $H\in \frak{a}$ such that $%
\alpha \left( H\right) =0$ if $\alpha \in \Theta $ and $\beta \left(
H\right) >0$ if $\beta \in \Pi ^{+}\setminus \langle {\Theta }\rangle $.
Then $\beta \in \Pi ^{+}\setminus \langle {\Theta }\rangle $ if and only if $%
\beta \left( \mathcal{C}_{\Theta }\right) >0$ (see Warner \cite{w}, Lemma
1.2.4.1). Also $w\in \mathcal{W}_{\Theta }$ if and only if $wH=H$, $H\in
\mathcal{C}_{\Theta }$ and the chambers meeting $\mathcal{C}_{\Theta }$ in
their closures are exactly $w\frak{a}^{+}$, $w\in \mathcal{W}_{\Theta }$.
Since roots do not change signs on chambers, it follows that $\alpha \in \Pi
^{+}\setminus \langle {\Theta }\rangle $ is $\geq 0$ on $\mathcal{W}_{\Theta
}\mathrm{cl}\frak{a}^{+}$.

Conversely we check that $w\mathrm{cl}\frak{a}^{+}$, $w\in \mathcal{W}%
_{\Theta }$, are the only closures of chambers where $\Pi ^{+}\setminus
\langle {\Theta }\rangle $ is $\geq 0$. In fact, take $w\notin \mathcal{W}%
_{\Theta }$ and $H\in \mathcal{C}_{\Theta }$. Then $wH\notin \mathcal{W}%
_{\Theta }\mathrm{cl}\frak{a}^{+}$, because otherwise $w$ would map $H$ to a
chamber containing $H$ in its closure and this can happen only if $w$ fixes $%
H$, that is only if $w\in \mathcal{W}_{\Theta }$. Now, the the set of roots
that are $>0$ on $wH$ is $w\left( \Pi ^{+}\setminus \langle {\Theta }\rangle
\right) $. By \cite{w}, Theorem 1.2.4.8 (and its proof) $w\left( \Pi
^{+}\setminus \langle {\Theta }\rangle \right) $ is not contained in $%
\langle {\Theta }\rangle \cup \Pi ^{+}$ if $w\notin \mathcal{W}_{\Theta }$.
Hence there exists $\gamma \in \Pi ^{-}\setminus \langle {\Theta }\rangle $
such that $\gamma \left( wH\right) >0$, that is, $\left( -\gamma \right)
\left( wH\right) <0$ so that $\Pi ^{+}\setminus \langle {\Theta }\rangle $
is not positive on $w\frak{a}^{+}$.

So far we have $\{H\in \frak{a}:\alpha \left( H\right) \geq 0,\alpha \in \Pi
^{+}\setminus \langle {\Theta }\rangle \}=\mathcal{W}_{\Theta }\mathrm{cl}%
\frak{a}^{+}$. But to take strict inequality $>0$ in the left hand side
amounts to take interior in the right hand side.

To prove the last statement we apply $w_2^{-1}$ to both cones and
reduce to the case where  $w_2 = 1$.  Now if $\mathrm{int}\left(
w_1 \mathcal{W}_{\Theta }\mathrm{cl}\frak{a}^{+}\right)$ meets
$\mathrm{int}\left( \mathcal{W}_{\Theta
}\mathrm{cl}\frak{a}^{+}\right)$ then, by the first part, they
meet at a regular element, say $w_1wH = w'H'$ where $w, w' \in
\mathcal{W}_{\Theta}$ and $H,H' \in \frak{a}^+$ are regular. Since
$\mathcal{W}$ acts simply on the Weyl chamber $\frak{a}^+$ it
follows that $H = H'$ so that $(w')^{-1}w_1w = 1$ since it fixes a
regular element.  This implies that $w_1 \in W_\Theta$, so that
both cones coincide.

\end{demonstracao}

By the above lemma we can restate Corollary \ref{corraizpositilambdamo} in
geometric terms as follows.

\begin{corolario}
\label{corincludecone}$\Lambda _{\mathrm{Mo}}(\mathcal{M}^{+})$ is
contained in the open convex cone $\mathrm{int}\left(
\mathcal{W}_{\phi }\mathrm{cl}\frak{a}^{+}\right) $.
\end{corolario}

\section{Morse spectra and block form}

In this section we combine the previous results to get the full
picture of the $\frak{a}$-Morse spectra of the several chain
components $\mathcal{M}\left( w\right) $ with $w$ running through
the Weyl group $\mathcal{W}$. As always we assume that the flow on
the base space $X$ is chain recurrent and let $\Theta =\Theta
({\phi}) $ be the parabolic type of the flow $\phi _{t}$.

It was proved before that the whole set of $\frak{a}$-Lyapunov
exponents is $\mathcal{W}$-invariant (see Corollary
\ref{corinvariance}). This result will be improved now by showing
that the Lyapunov spectra of the chain components are permuted to
each other by the Weyl group. To this end we recall the $Z_{\Theta
}$-block reduction $Q_{\phi}$ of Section \ref{secpartydyn} as well
as the $K_{\phi }$-reduction $R_{\phi }$. If $b_{0}$ stands for
the origin in $\Bbb{F}$ then by Proposition
\ref{propchaincompdynreduct} we have
\[
{\mathcal{M}}(w)=Q_\phi\cdot wb_{0}=R_{\Theta }\cdot wb_{0}.
\]

\begin{lema}
Let $r\in R_{\Theta }$ be regular. Then it has asymptotic ray $D=\mathrm{Ad}%
\left( u\right) H^{+}$ where $u\in K_{\Theta }$ and $H^{+}\in \mathrm{cl}%
\frak{a}^{+}$ is the polar exponent of $r$. Moreover,  $\Theta
\left( H^{+}\right) \subset \Theta $.
\end{lema}

\begin{demonstracao}
Since $\phi _{k}(r)\in Q_{\phi }$ we decompose $\phi _{k}\left(
r\right) =r_{k}\cdot s_{k}\in R_{\Theta }\cdot \left( S\cap
Z_{\Theta }\right) $, then $s_{k}=\mathsf{S}(\phi _{k}(r))\in
Z_{\Theta }$ is regular with asymptotic ray $D\in \frak{s}$. By
Proposition \ref{prop:kaimanovich-em-Z} we have
$D=\mathrm{Ad}\left( u\right) H$ where $u\in K_{\Theta }$ and
$H\in \frak{a}$ is such that $\alpha (H)\geq 0$ for all $\alpha
\in \Theta $.

Noting that $\mathsf{S}(\phi _{k}(r\cdot u))=u^{-1}s_{k}u$ is a
regular sequence in $G$ with asymptotic ray $H=\mathrm{Ad}\left(
u^{-1}\right) D$ it follows that $r\cdot u\in R_{\Theta }$ is a
regular point of the flow with
asymptotic ray $H\in \frak{a}$. Therefore we have the Lyapunov exponent $%
\lambda (ru\cdot b_{0})=H$ (see Proposition
\ref{propkaimanovichflow} and
Proposition \ref{propkaimanovich} (2)). Since $ru\cdot b_{0}\in \mathcal{M}%
^{+}$ it follows that $H\in \Lambda
_{\mathrm{Ly}}(\mathcal{M}^{+})$ so that, applying the estimates
of last section, we conclude that $\alpha (H)>0$ for all $\alpha
\in \Pi ^{+}\setminus \langle {\Theta }\rangle $ (see Corollary
\ref{corparabtyplypexp}). It follows  that $\Theta \left( H\right)
\subset \Theta $. Also,  $\alpha (H)\geq 0$ for all $\alpha \in
\Pi ^{+}$ so that $H\in \mathrm{cl}\frak{a}^{+}$. This implies
that $H$ is the polar exponent of $r$ which is denoted by
$H^{+}=H$ in the statement of the result.
\end{demonstracao}

\begin{proposicao}
Let $r\in R_{\Theta }$ be regular with asymptotic ray $D\in
\frak{s}$ and polar exponent $H^{+}\in \mathrm{cl}\frak{a}^{+}$.
Then for $x=\pi (r)$ we have
\[
\bigcup_{s\in \mathcal{W}_{\Theta _{\phi }}}r\cdot
\mathrm{fix}(D,sw)\subset \mathcal{M}(w)_{x},
\]
\[
\Lambda
_{\mathrm{Ly}}(\mathcal{M}(w)_{x})=w^{-1}\mathcal{W}_{\Theta
_{\phi }}H^{+}.
\]
\end{proposicao}

\begin{demonstracao}
By the previous lemma we have $D=\mathrm{Ad}\left( u\right) H^{+}$ with $%
u\in K_{\Theta }$, $H^{+}\in \mathrm{cl}\frak{a}^{+}$ and $\Theta
\left( H^{+}\right) \subset \Theta $. Let $H_{\phi }$ be as in
Theorem \ref {teo:tipo-parabol-flow} so that $\Theta
(H^{+})\subset \Theta =\Theta
(H_{\phi })$. Hence for any $w\in \mathcal{W}$, $s\in \mathcal{W}_{\Theta }$%
, we have $\mathrm{fix}\left( H^{+},sw\right) \subset
\mathrm{fix}\left( H_{\phi },sw\right) =\mathrm{fix}\left( H_{\phi
},w\right) $, so that
\[
r\cdot \mathrm{fix}(D,sw)=ru\cdot \mathrm{fix}(H^{+},sw)\subset
ru\cdot \mathrm{fix}(H_{\phi },w)=\mathcal{M}(w).
\]
By  Proposition \ref{propkaimanovich} (3)  the above inclusion
implies that
\begin{equation}
w^{-1}\mathcal{W}_{\Theta }H^{+}\subset \Lambda _{\mathrm{Ly}}(\mathcal{M}%
(w)_{x}).  \label{eq4:permutacao-espectro}
\end{equation}
Now let $\xi \in \mathcal{M}(w)_{x}$ and write $\xi =r\cdot
lwb_{0}$, $l\in K_{\Theta }$. By  Proposition
\ref{propkaimanovich} (3) we have $\lambda (\xi )=s^{-1}H^{+}$
where $s\in \mathcal{W}$ is such that $lwb_{0}\in
\mathrm{st}(D,s)$. Since $u\in K_{\Theta }$ we have
\[
lwb_{0}\in \mathrm{st}(X,s)=uP_{\Theta \left( H^{+}\right)
}^{-}sb_{0}\subset P_{\Theta }^{-}sb_{0},
\]
and since $l\in K_{\Theta }$ one has that $wb_{0}\in P_{\Theta
}^{-}sb_{0}$.
By the Bruhat decomposition it follows that $\mathcal{W}_{\Theta }w=\mathcal{%
W}_{\Theta }s$ so that $s^{-1}\in w^{-1}\mathcal{W}_{\Theta }$.
Hence
\[
\lambda (\xi )=s^{-1}H^{+}\in w^{-1}\mathcal{W}_{\Theta }H^{+},
\]
which implies the reverse inclusion in
(\ref{eq4:permutacao-espectro}).
\end{demonstracao}

Now we are ready to state the permutation of the Lyapunov and
Morse spectra of the Morse components under the Weyl group.

\begin{teorema}
\label{teo:permutacao-espectro}For every $w\in \mathcal{W}$ we
have
\[
\Lambda _{\mathrm{Ly}}\left( {\mathcal{M}}\left( w\right) \right)
=w^{-1}\Lambda _{\mathrm{Ly}}({\mathcal{M}}^{+})\qquad
\mathrm{and}\qquad \Lambda _{\mathrm{Mo}}\left(
{\mathcal{M}}\left( w\right) \right) =w^{-1}\Lambda
_{\mathrm{Mo}}({\mathcal{M}}^{+}).
\]
\end{teorema}

\begin{demonstracao}
The statement for the Lyapunov spectra implies that for the Morse
spectra because $\Lambda _{\mathrm{Mo}}\left( {\mathcal{M}}\left(
w\right) \right) $ is the convex closure of $\Lambda
_{\mathrm{Ly}}\left( {\mathcal{M}}\left( w\right) \right) $.

For the Lyapunov spectra we recall that  any Lyapunov exponent
$\lambda $ is $\lambda =\lambda \left( r\cdot b\right) =\lambda
\left( s_{k},b\right) $
for  a regular point $r\in R$ and $b\in \Bbb{F}$  where $s_{k}=\mathsf{S}%
\left( \phi _{k}\left( r\right) \right) $ (see Proposition \ref
{propkaimanovichflow} and Corollary
\ref{corol:espectro-pontos-reg}). Hence,
\begin{eqnarray*}
\Lambda _{\mathrm{Ly}}(\mathcal{M}(w))
&=&\bigcup_{x\,\mathrm{regular}\text{
}}\Lambda _{\mathrm{Ly}}(\mathcal{M}(w)_{x})=\bigcup_{x\,\mathrm{reg}}w^{-1}%
\mathcal{W}_{\Theta }\lambda ^{+}(x)= \\
&=&w^{-1}\left( \bigcup_{x\,\mathrm{reg}}\mathcal{W}_{\Theta
}\lambda ^{+}(x)\right) =w^{-1}\Lambda
_{\mathrm{Ly}}(\mathcal{M}^{+}),
\end{eqnarray*}
which proves the theorem.
\end{demonstracao}

\begin{remar}
The above theorem says that the several Morse spectra are
homeomorphic to each other. This is in consonance with the fact
that the intersection of the different chain components with the
fibers are  homeomorphic as well.  In fact, a such intersection
can be identified to the orbit  $K_{H_{\phi }}wb_{0}$, which in
turn identifies to the coset space $K_{H_{\phi }}/M$.
\end{remar}

Combining this theorem with the inclusion in chambers stated in Corollary
\ref{corincludecone} we get the following localization result for the Morse
(and Lyapunov) spectra.

\begin{corolario}
For each $w\in \mathcal{W}$ the Lyapunov and Morse spectra $\Lambda _{%
\mathrm{Ly}}\left( {\mathcal{M}}\left( w\right) \right) $ and $\Lambda _{%
\mathrm{Mo}}\left( {\mathcal{M}}\left( w\right) \right) $ are contained in
the open cone
\[
\mathrm{int}\left( w^{-1}\mathcal{W}_{\Theta_{phi}
}\mathrm{cl}\frak{a}^{+}\right) .
\]
It follows that $\Lambda _{\mathrm{Mo}}\left( {\mathcal{M}}\left(
w_{1}\right) \right) \cap \Lambda _{\mathrm{Mo}}\left( {\mathcal{M}}\left(
w_{2}\right) \right) =\emptyset $ if $\mathcal{M}\left( w_{1}\right) \neq
\mathcal{M}\left( w_{2}\right) $, that is, if $\mathcal{W}_{\phi }w_{1}\neq
\mathcal{W}_{\phi }w_{2}$.
\end{corolario}

\begin{demonstracao}
In fact, $\mathrm{int}\left( w^{-1}\mathcal{W}_{\Theta ({\phi}) }
\mathrm{cl}\frak{a}^{+}\right) =w^{-1}$ $\mathrm{int}\left(
\mathcal{W}_{\Theta ({\phi}) }\mathrm{cl}\frak{a}^{+}\right) $
and the right hand side contains $\Lambda _{\mathrm{Mo}}({\mathcal{M}}^{+})$%
. The last statement follows from the fact that the open cones $\mathrm{int}%
\left( w_{1}^{-1}\mathcal{W}_{\Theta ({\phi}) }
\mathrm{cl}\frak{a}^{+}\right) $ and $\mathrm{int}\left( w_{2}^{-1}%
\mathcal{W}_{\Theta ({\phi}) }\mathrm{cl}\frak{a}%
^{+}\right) $ are either equal or disjoint.
\end{demonstracao}

In particular, the equalities in Theorem \ref{teo:permutacao-espectro} say
that if $w\in \mathcal{W}_{\phi }$, that is, if $\mathcal{M}\left( w\right) =%
\mathcal{M}\left( 1\right) =\mathcal{M}^{+}$, then $w\Lambda _{\mathrm{Ly}}({%
\mathcal{M}}^{+})=\Lambda _{\mathrm{Ly}}({\mathcal{M}}^{+})$ and $w\Lambda _{%
\mathrm{Mo}}({\mathcal{M}}^{+})=\Lambda _{\mathrm{Mo}}({\mathcal{M}}^{+})$.
This invariance is easily carried over to the other spectra, by taking a
conjugate of $\mathcal{W}_{\phi }$. Namely if $s\in w^{-1}\mathcal{W}_{\phi
}w$ then $s\Lambda _{\mathrm{Ly}}\left( {\mathcal{M}}\left( w\right) \right)
=\Lambda _{\mathrm{Ly}}\left( {\mathcal{M}}\left( w\right) \right) $ and $%
s\Lambda _{\mathrm{Mo}}\left( {\mathcal{M}}\left( w\right) \right)
=\Lambda _{\mathrm{Mo}}\left( {\mathcal{M}}\left( w\right) \right)
$, as follows directly from the theorem. Actually, we have the
following more precise result.

\begin{corolario}
For each $w\in \mathcal{W}$ we have $w^{-1}\mathcal{W}_{\phi }w=\{s\in
\mathcal{W}:s\Lambda _{\mathrm{Mo}}\left( {\mathcal{M}}\left( w\right)
\right) =\Lambda _{\mathrm{Mo}}\left( {\mathcal{M}}\left( w\right) \right) \}
$. The same statement holds with $\Lambda _{\mathrm{Ly}}$ instead of $%
\Lambda _{\mathrm{Mo}}$.
\end{corolario}

\begin{demonstracao}
In fact, $s\in \mathcal{W}$ fixes $\Lambda _{\mathrm{Mo}}\left( {\mathcal{M}}%
\left( w\right) \right) $ if and only if $wsw^{-1}$ fixes $\Lambda _{\mathrm{%
Mo}}\left( {\mathcal{M}}^{+}\right) $. By the above theorem $\Lambda _{%
\mathrm{Mo}}\left( {\mathcal{M}}^{+}\right) $ is invariant under $\mathcal{W}%
_{\phi }$. Conversely, if $\sigma \Lambda _{\mathrm{Mo}}\left( {\mathcal{M}}%
^{+}\right) =\Lambda _{\mathrm{Mo}}\left( {\mathcal{M}}^{+}\right) $ then $%
\Lambda _{\mathrm{Mo}}\left( {\mathcal{M}}\left( \sigma ^{-1}\right) \right)
=\Lambda _{\mathrm{Mo}}\left( {\mathcal{M}}^{+}\right) $ hence $\sigma \in
\mathcal{W}_{\phi }$.
\end{demonstracao}

In view of this corollary the parabolic subgroup $\mathcal{W}_{\phi }$ can
be recovered from the Morse spectra of $\phi $ (more specifically from the
Morse spectrum $\Lambda _{\mathrm{Mo}}\left( {\mathcal{M}}^{+}\right) $). In
other words the parabolic type of $\phi $ is completely determined by its
spectra.

Another way of expressing the relationship between the exponents
and the parabolic type is by comparing $\Theta (\phi)$ with the
set of simple roots annihilating some Morse exponents, as stated
in the next consequence of Theorem \ref{teo:permutacao-espectro}.

\begin{corolario}
For a Morse exponent $\lambda \in \Lambda _{\mathrm{Mo}}\left( \mathcal{M}%
^{+}\right) $ write $\Theta (\lambda )=\{\alpha \in \Sigma :\alpha
\left( \lambda \right) =0\}$. Then $\Theta (\lambda )\subset
\Theta (\phi) $. Also there exists $\lambda _{0}\in \Lambda
_{\mathrm{Mo}}\left( \mathcal{M}^{+}\right) $ such that $\Theta
(\lambda )\subset \Theta (\phi) $.
\end{corolario}

\begin{demonstracao}
The inclusion $\Theta (\lambda )\subset \Theta (\phi) $ is just a
restatement of Corollary \ref{corraizpositilambdamo}. On the other
hand
the $\mathcal{W}_{\phi }$-invariance of convex set $\Lambda _{\mathrm{Mo}%
}\left( \mathcal{M}^{+}\right) $ implies this set contains a point
$\lambda _{0}$ fixed by $\mathcal{W}_{\phi }$. This fixed point
satisfies $\Theta (\lambda _{0})\supset \Theta (\phi) $, hence the
equality.
\end{demonstracao}

In other words the block form of $\phi $ is given by the elements
in $\Lambda _{\mathrm{Mo}}\left( \mathcal{M}^{+}\right) $ with the
smallest possible degree of regularity. Hence, again the block
form of $\phi $ can be read off from the its Morse spectrum.

As an immediate consequence we mention that the parabolic type $\Theta
=\Sigma $ (or block form $H_{\phi }=0$) corresponds to chain transitivity of
the flow on the flag bundles. Hence we have the following criterion for
chain transitivity on the flag bundles.

\begin{proposicao}
\label{aplic:espectro-do-atrator}The flow $\phi _{t}$ is chain transitive in
some flag bundle if and only if $0\in \Lambda _{\mathrm{Mo}}(\mathcal{M}^{+})
$. In this case $\phi _{t}$ is chain transitive (and chain recurrent) on
every flag bundle.
\end{proposicao}

\begin{exemplo}
If the base $X$ is a point then our  flow is just the  iteration
of a single element $g\in G$  (in the discrete time case) or  the
action of a one-parameter group $\exp (tX)$ in $G$, $X\in
\frak{g}$  (continuous time). For the discrete-time case generated
by $g\in G$, let $g=uhn$ be its Jordan-Schur decomposition (see
\cite{Helgason}, Chapter IX), where $u$ is
elliptic, $h$ is hyperbolic, and $u$ unipotent. Choose a Weyl chamber $A^{+}$%
of $G$ such that $h\in \mathrm{cl}A^{+}$and put $h=\exp (H^{+})$,
$H^{+}\in \mathrm{cl}\frak{a}^{+}$. Then it can be shown that  the
block form of this flow is  the block form of the hyperbolic part
$H^{+}$ (see In \cite {psmsec-jordan}). Also  the finest Morse
decomposition in $\Bbb{F}$ is given
by the fixed point set  of $h$, that is, $\mathcal{M}(w)=\mathrm{fix}%
(H^{+},w)$. To compute the vector exponents of this flow we fix a
Cartan and Iwasawa decomposition of $G$ compatible with the
chamber $A^{+}$. The polar exponent of $g^{k}$ is precisely
$H^{+}$ (cf.   \cite{tamtiny-svd}, Theorem 2.2) . So that by our
results it follows that
\[
\Lambda _{\mathrm{Ly}}(\mathcal{M}(w))=\Lambda _{\mathrm{Mo}}(\mathcal{M}%
(w))=w^{-1}H^{+}.
\]
\end{exemplo}

\section{Representations and vector bundles\label{secvecbundl}}

In this final section we establish some relationships between the
$\frak{a} $-Lyapunov and Morse spectra and the exponents of a
linear flow on a vector bundle.

Let $\mathcal{V}\rightarrow X$ be a finite dimensional real vector
bundle with compact Hausdorff base space and $\phi _{t}$ a linear
flow on $\mathcal{V}$, assumed to be chain transitive on the base
space. The bundle $Q=B\mathcal{V}\rightarrow X$ of frames of
$\mathcal{V}$ is a principal bundle with structural group
$\mathrm{Gl}\left( m,\Bbb{R}\right) $,
$m=\mathrm{rank}\,\mathcal{V}$. The linear flow $\phi _{t}$ maps
frames into frames, hence it lifts to a right invariant flow (also
denoted by $\phi _{t}$) on $B\mathcal{V}$.

Fix in $\mathcal{V}$ a Riemannian metric $\langle \cdot ,\cdot
\rangle $ and let $O\mathcal{V}$ be the $\mathrm{O}\left( n\right)
$-subbundle of orthonormal frames of $B\mathcal{V}$. Then
$O\mathcal{V} \cdot \overline{A}N$ and $O\mathcal{V} \cdot S$ are,
respectively, the Iwasawa and Cartan decompositions of
$B\mathcal{V}$, where $\overline{A}$ is the subgroup of diagonal
matrices (with positive entries), $N$ is the subgroup of upper
triangular unipotent matrices and $S$ the set of positive definite
matrices of $\mathrm{Gl}\left( m,\Bbb{R}\right)$.

The spectra are defined by the cocycle $\mathsf{a}\left( t, r\right)
=\log a_{t}\in \overline{\frak{a}}$, where $r \in O\mathcal{V}$,
$\overline{\frak{a}}$ is the subspace of diagonal matrices and
$a_{t}$ is the $\overline{A}$-component in the Iwasawa decomposition
of $\phi_t(r) = r_t \cdot a_t n_t \in O\mathcal{V} \cdot
\overline{A}N$. This cocycle descends to the flag bundle
$\Bbb{F}\mathcal{V}$ of the flags of subspaces of $\mathcal{V}$.

Write $\overline{\frak{a}}=\Bbb{R}\cdot \mathrm{id}\oplus \frak{a}$,
where $\frak{a}$ is the subspace of zero trace diagonal matrices.
Then the component of the spectra in the direction of the scalar
matrices $\Bbb{R}\cdot \mathrm{id}$ is the same for any Morse
component in $\Bbb{F}\mathcal{V}$ (see Section \ref{sec-reductive}),
while the $\frak{a}$-component is
invariant under the Weyl group (permutation group of $m$ letters).

Now, let $\Bbb{P}\mathcal{V}$ be the projective bundle of $\mathcal{V}$. The
spectra of the flow on $\mathcal{V}$ is defined by the cocycle $\mathsf{a}%
_{\left| \cdot \right| }\left( t,v\right) $ over $\Bbb{P}\mathcal{V}$ given
by
\[
\mathsf{a}_{\left| \cdot \right| }\left( t,v\right) =\log \frac{\left| \phi
_{t}v\right| }{\left| v\right| }\qquad v\in \mathcal{V}\setminus \{0\}.
\]
To get the relation between the two cocycles take $v\in \mathcal{V}$ with $%
\left| v\right| =1$ and write $v=r\cdot e_{1}$, $r\in O\mathcal{V}$ and $%
e_{1}\in \Bbb{R}^{m}$ the first basic vector. Then $\phi _{t}\left(
v\right) =\phi _{t}\left( r\right) \cdot e_{1}$. Take the Iwasawa
decomposition $\phi _{t}\left( r\right) =r_{t}\cdot a_{t}n_{n}\in
\mathrm{O}\mathcal{V} \cdot AN$. Then $\left| \phi _{t}v\right|
=\left| a_{t}n_{t}e_{1}\right| =\left| a_{t}e_{1}\right| $. So that
\[
\mathsf{a}_{\left| \cdot \right| }\left( t,v\right) =\lambda _{1}\left( \log
a_{t}\right)
\]
where $\lambda _{1}:\frak{a}\rightarrow \Bbb{R}$ is the first eigenvalue $%
\lambda _{1}\left( \mathrm{diag}\{a_{1},\ldots ,a_{m}\}\right) =a_{1}$.
Hence
\[
\mathsf{a}_{\left| \cdot \right| }\left( t,v\right) =\lambda
_{1}\left( \mathsf{a}\left( t, \xi \right) \right)
\]
for any $\xi\in \Bbb{F}\mathcal{V}$ projecting into $[v] \in
\Bbb{P}\mathcal{V}$. On the other hand a chain component on the
projective bundle is the projection of chain components on the flag
bundle. This implies the following expressions for the spectra of
the linear flow on the vector bundle.

\begin{proposicao}
Let $\mathcal{M}\subset \Bbb{P}\mathcal{V}$ be a chain component and denote
by $\Lambda _{\mathrm{Mo}}^{\left| \cdot \right| }\left( \mathcal{M}\right) $
and $\Lambda _{\mathrm{Ly}}^{\left| \cdot \right| }\left( \mathcal{M}\right)
$ its spectra. Then $\Lambda _{\mathrm{Mo}}^{\left| \cdot \right| }\left(
\mathcal{M}\right) =\lambda _{1}\left( \Lambda _{\mathrm{Mo}}\left(
\overline{\mathcal{M}}\right) \right) $ and $\Lambda _{\mathrm{Ly}}^{\left|
\cdot \right| }\left( \mathcal{M}\right) =\lambda _{1}\left( \Lambda _{%
\mathrm{Ly}}\left( \overline{\mathcal{M}}\right) \right) $, where $\overline{%
\mathcal{M}}$ is any chain component on $\Bbb{F}\mathcal{V}$ projecting onto
$\mathcal{M}$.
\end{proposicao}

In particular if the block form of $\phi _{t}$ is a regular matrix
then there are $m!$, $m=\mathrm{rank}\mathcal{V}$, chain components
on the maximal flag bundle $\Bbb{F}\mathcal{V}$. Each one meets the
fibers in singletons. By the previous results (see Corollary
\ref{corincludecone} and Theorem \ref{teo:permutacao-espectro}) each
Morse spectra of a chain component $\mathcal{M}\left( w\right) $ is
entirely contained in the Weyl chamber $w\frak{a}^{+}$, and so
contains only regular elements. The chain components on
$\Bbb{F}\mathcal{V}$ project down to $m$ components on the
projective bundle $\Bbb{P}\mathcal{V}$, which are ordered linearly.
The vector bundle Morse spectrum $\Lambda _{\mathrm{Mo}}^{\left|
\cdot \right| }$ of the attractor component is just the set of
highest eigenvalues of the matrices in $\Lambda _{\mathrm{Mo}}\left(
\mathcal{M}^{+}\right) $. For the other
components one must apply a permutation to the elements of $\Lambda _{\mathrm{%
Mo}}\left( \mathcal{M}^{+}\right) $ and then take the highest eigenvector.
This amounts to take the successive eigenvectors of the elements of $\Lambda
_{\mathrm{Mo}}\left( \mathcal{M}^{+}\right) $. Since the elements of $%
\Lambda _{\mathrm{Mo}}\left( \mathcal{M}^{+}\right) $ (and hence of $\Lambda
_{\mathrm{Ly}}\left( \mathcal{M}^{+}\right) $) are regular matrices this
yields the simplicity of the Lyapunov spectrum of a nominal trajectory on
the base space.

We note that the  above proposition can be easily extended to
obtain the spectra on the Grassmann bundles
$\mathrm{Gr}_k\mathcal{V}$ from the vector spectra on the flag
bundle $\Bbb{F}\mathcal{V}$: One needs only to replace
$\Bbb{P}\mathcal{V}$ by $\mathrm{Gr}_k\mathcal{V}$ and the
functional $\lambda_1$ by the
functional $\lambda_k$ which gives the sum of the $k$-first eigenvalues %
$\lambda _{k}\left( \mathrm{diag}\{a_{1},\ldots ,a_{m}\}\right)
=a_{1}+\cdots+a_k$.

In the general case the idea is to start with a principal bundle
$Q\rightarrow X$ with structural group $G$ and take a linear
representation $\rho :G\rightarrow \mathrm{Gl}\left( V\right) $ on a
vector space $V$. This representation yields the associated bundle
$\mathcal{V}=Q\times _{\rho }V\rightarrow X$ which turns out to be a
vector bundle. If $\phi _{t}$ is a right invariant flow on $Q$ then
the induced flow on $\mathcal{V}$ is linear.

For instance if $G$ is semi-simple and connected then an irreducible
representation $\rho _{\mu }$ is given by its highest weight $\mu \in \frak{a%
}^{*}$. Then there is a $K$ invariant inner product on $\mathcal{V}$ such
that the norm cocycle $\mathsf{a}_{\left| \cdot \right| }\left( t,v\right)
=\log \frac{\left| \phi _{t}v\right| }{\left| v\right| }$ is given by
\[
\mathsf{a}_{\left| \cdot \right| }\left( t,v\right) =\mu \left( \log
a_{t}\right)
\]
with $a_{t}$ the $A$-component in the Iwasawa decomposition
of $\phi_t(r)$ where $v = r \cdot Y$ for a highest weight vector
$Y \in V$.

Also, for $h\in \mathrm{cl}A^{+}$, $\log \Vert \rho _{\mu }(h)\Vert
=\langle \mu ,\log h\rangle $ (cf. \cite{gjt}, Lemma 4.22). Thus the
vector bundle spectra of the linear flow on $\mathcal{V}$ can be
recovered from the intrinsic spectra on the principal bundle. We do
not go here into the details, since it requires a discussion of the
chain components on the projective bundle of $\mathcal{V}$, which is
not yet available.


\begin{thebibliography}{99}
\bibitem{aco}  Arnold, L., N. D. Cong and V.I. Oseledets: \textit{Jordan
normal form for linear coccycles}. Random Oper. Stochastic
Equations \textbf{7} (1999), 303-358.

\bibitem{Ayala1} Ayala,V.,  F. Colonius and W. Kliemann, \textit{Dynamical
characterization of the Lyapunov form of matrices}, Linear Algebra
and Its Applications 420(2005), 272 - 290.

\bibitem{Ayala2} Ayala, V.,  F. Colonius and W. Kliemann: \textit{On topological
equivalence of linear flows with applications to bilinear control
systems}, submitted.

\bibitem{smbflow}  Braga Barros, C.J. and L.A.B. San Martin: \textit{Chain
transitive sets for flows on flag bundles}. Forum Math., vol. 19 (2007),
19-60.

\bibitem{ck1}  Colonius, F. and Kliemann, W.: \textit{The Morse spectrum of
linear flows on vector bundles}. Trans. AMS \textbf{348} (1996), 4355-4388.

\bibitem{ck}  Colonius, F. and Kliemann, W.: The dynamics of control.
Birkh\"{a}user, Boston (2000).

\bibitem{cfj}  Colonius, F., R. Fabbri and R.A. Johnson: \textit{Chain
Recurrence, Growth Rates and Ergodic Limits}, Ergodic Theory and Dynamical
Systems, to appear.

\bibitem{CFJS} Colonius, F., R.  Fabbri,  R.A. Johnson and M. Spadini:
\textit{Bifurcation phenomena in control flows}, Topological
Methods in Nonlinear Analysis, (to appear).

\bibitem{c}  Conley C.: \textit{Isolated invariant sets and the Morse
index}.  CBMS Regional Conf. Ser. in Math., \textbf{38}, American
Mathematical Society, (1978).

\bibitem{c1}  Conley C.: \textit{The gradient structure of a flow: I}.
Ergodic Theory Dynam. Systems, \textbf{8} (1988), 11-26.

\bibitem{dkv}  Duistermat, J.J., J.A.C. Kolk and V.S. Varadarajan: \textit{%
Functions, flows and oscilatory integral on flag manifolds}. Compositio
Math. 49, 309-398, (1983).

\bibitem{fer}  Feres, R.: Dynamical systems and semisimple groups. An
introduction. Cambridge University Press (1998).

\bibitem{gjt}  Y. Guivarc'h, L. Ji and J. C. Taylor: Compactifications of
Symmetric Spaces. Springer-Verlag (1998).

\bibitem{Helgason} S. Helgason: Differential Geometry, Lie Groups and
Symmetric Spaces, Academic Press, (1978).

\bibitem{jps}  Johnson, R.A., K.J. Palmer and G.R. Sell: \textit{Ergodic
properties of linear dynamical systems}, SIAM J. Math. Anal.
\textbf{18}, (1987), 1-33.

\bibitem{ka}  Kaimanovich, V.A.: \textit{Lyapunov exponents, symmetric
spaces, and a multiplicative ergodic theorem for semisimple Lie groups}. J.
Soviet Math. \textbf{47} (1989), 2387-2398.

\bibitem{kapp}  Knapp, A.W.: Lie Groups Beyond an Introduction. Birkh\"{a}user, (2002).

\bibitem{kn}  Kobayashi, S. and K. Nomizu: Foundations of Differential
Geometry vol.I, InterScience Publishers, (1963).

\bibitem{Gabriele} Link, G.: \textit{Limit Sets of Discrete Groups
acting on Symmetric Spaces}, PhD thesis, Fakult\"{a}t f\"{u}r
Mathematik der Universit\"{a}t Karlsruhe, (2002).

\bibitem{patr}  Patr\~{a}o, M.: \textit{Morse decomposition of semiflows on
topological spaces}. J. Dyn. Diff. Eq., \textbf{19} (2007), 181-198.

\bibitem{msm1}  Patr\~{a}o, M. and L.A.B. San Martin: \textit{Semiflows on
Topological Spaces: Chain transitivity and Shadowing Semigroups}. J. of
Dynamics and Diff. Eq. \textbf{19} (2007) 155-180.

\bibitem{msm}  Patr\~{a}o, M. and L.A.B. San Martin: \textit{Chain
recurrence of flows and semiflows on fiber bundles}. Discrete Contin. Dynam.
Systems A, \textbf{17} (2007), 113-139.

\bibitem{psmsec}  Patr\~{a}o, M., L.A.B. San Martin and L. Seco, L.: \textit{Stable
manifolds and Conley index for flows in flag bundles}. (to appear).

\bibitem{psmsec-jordan}  Patr\~{a}o, M., L.A.B. San Martin and L. Seco:
\textit{Parabolic type and the Jordan decomposition}. (to appear)

\bibitem{salzh}  Salamon, D. and E. Zehnder: \textit{Flows on vector bundles
and hyperbolic sets}. Trans. AMS \textbf{306} (1988), 623-659.

\bibitem{sm}  San Martin, L.A.B.: \textit{Invariant control sets on flag
manifolds}. Math. Control Signals Systems 6 no. 1, 41-61, (1993).

\bibitem{smt}  San Martin, L.A.B. and P.A. Tonelli: \textit{Semigroup
actions on homogeneous spaces. } Semigroup Forum, \textbf{50} (1995), 59-88.

\bibitem{sl}  Seco, L.: \textit{Morse exponents of vector valued cocycles}.
(to appear)

\bibitem{tamtiny-svd} Tam, T.Y. and  H. Huang: \textit{An
extension of Yamamoto's theorem on the eigenvalues and singular
values of a matrix}, Journal of Math. Soc. Japan no.58, 1197-1202,
(2006).

\bibitem{w}  Warner, G.: Harmonic analysis on semi-simple Lie groups I.
Springer-Verlag (1972).

\bibitem{z}  Zimmer, R.: Ergodic theory and semisimple groups. Monographs in
Mathematics. Birkh\"{a}user (1984).
\end{thebibliography}
\end{document}